\documentclass[12pt]{amsart}
\usepackage[utf8]{inputenc}
\usepackage[foot]{amsaddr}
\usepackage{amsmath,amsfonts,amsthm,amssymb,verbatim,etoolbox,color}
\usepackage{enumitem} 
\usepackage{epigraph}
\usepackage{mathtools}
\usepackage{appendix}
\usepackage[normalem]{ulem}
\usepackage{xfrac}
\usepackage{flexisym}
\usepackage{bbm}
\usepackage{comment}
\usepackage{tikz}
\usetikzlibrary{patterns}

\usepackage[margin=1in]{geometry}
\usepackage[
    backend=biber,
    maxnames=8,
    maxalphanames=8,
    style=alphabetic,
    backref,
  ]{biblatex}
\usepackage{scalerel}

\usepackage{hyperref}
\hypersetup{
    colorlinks=true,
    linkcolor=red,
    filecolor=magenta,      
    urlcolor=magenta,
    pdftitle={},
    pdfauthor = {Ihor Pylaiev},
    linktocpage = true,
    citecolor = blue,
    }

\usepackage{theoremref}

\newtheorem{thm}{Theorem}
\newtheorem{lemma}[thm]{Lemma}
\newtheorem{lem}[thm]{Lemma}

\newtheorem{corollary}[thm]{Corollary}

\newtheorem{prop}[thm]{Proposition}
\theoremstyle{definition}
\newtheorem{defn}[thm]{Definition}
\newtheorem{remark}[thm]{Remark}
\newtheorem{claim}[thm]{Claim}

\numberwithin{thm}{section}

\newcommand{\N}{\mathbb{N}}

\newcommand{\Z}{\mathbb{Z}}
\newcommand{\R}{\mathbb{R}}
\renewcommand{\P}{\mathbb{P}}

\newcommand{\Q}{\mathbb{Q}}
\newcommand{\C}{\mathbb{C}}
\newcommand{\T}{\mathbb{T}}
\newcommand{\EE}{\mathbb{E}}
\DeclareMathOperator{\rad}{rad}
\newcommand{\doublefactorial}{!\!!}

\newcommand{\bra}[1]{\left(#1\right)}
\newcommand{\sqbra}[1]{\left[#1\right]}

\newcommand{\blue}[1]{#1}

\newcommand{\ra}{\rightarrow}

\newcommand{\conj}[1]{\overline{#1}}
\newcommand{\abs}[1]{\left\lvert {#1} \right \rvert}
\newcommand{\aabs}[1]{\left\lVert {#1} \right \rVert}

\newcommand{\var}{\text{Var}}
\newcommand{\cov}{\text{Cov}}

\newcommand{\Span}{\text{Span}}

\newcommand{\ep}{\varepsilon}

\newcommand{\HBV}{\text{HBV}}

\newcommand{\Pp}{P^+}
\newcommand{\Ppsf}{P_{\text{sf}}^+}

\title{Random Chowla and Goldbach Conjectures}
\author{Ihor Pylaiev}
\email{ipylaiev@mit.edu}
\address{Department of Mathematics, Massachusetts Institute of Technology}
\thanks{\blue{I have used Generative AI tools for polishing the writing in this paper and literature search, but I assume full responsibility for the originality, mathematical correctness, and final content of this work.}}
\date{\today}

\begin{document}

\begin{abstract}
We prove that for a Rademacher random multiplicative function $(f(n))_{n\in\N}$ and any non-linear polynomial $P\in\Z[x]$ admitting a positive constant $\kappa_P$ such that
\[\lim_{N\rightarrow\infty}
\frac{\#\{n\in\sqbra{N}:P(n)\text{ is square-free}\}}{N}
=\kappa_P,\]
we have 
\[
\frac{1}{\sqrt{\kappa_PN}}\sum_{n\leqslant N}f(P(n))
\xrightarrow[]{d}\mathcal{N}(0,1),
\]
as $N\rightarrow\infty $, where convergence is in distribution to a standard Gaussian. Using a direct combinatorial argument, we establish paucity estimates for non-trivial square products of square-free polynomial values at every fixed moment order. These estimates yield the Gaussian limit by the method of moments.

We also study partial sums of the form $\sum_{n<N} F(f(n),f(N-n))$ for a suitable class of bivariate functions $F$ when $f$ is a Steinhaus random multiplicative function, and prove analogous results in special cases in the Rademacher setting. These results can be viewed as random analogues of the Goldbach-type problems considered by Mangerel.
\end{abstract}

\maketitle
\tableofcontents

\section{Introduction}
The distribution of multiplicative functions, particularly the Möbius function $\mu(n)$ plays a crucial role in the theory of prime numbers. In particular, many problems in Analytic Number Theory can be rephrased to state that $\mu(n)$ behaves pseudo-randomly. For example,
\[\sum_{n=1}^N\mu(n)=o(N)\]
is equivalent to the Prime Number Theorem, and the well-known Riemann Hypothesis is equivalent to 
\[\sum_{n=1}^N \mu(n)=O_\epsilon(N^{\frac{1}{2}+\epsilon})\text{ for any $\epsilon>0$},\]
which is almost the same level of cancellation as we would expect from independent uniform $\pm 1$ random variables.
\blue{\begin{remark}
    Throughout the paper, I will write $O_\epsilon(\cdot)$ or $\ll_\epsilon$ to denote that the bound holds for any $\epsilon>0$, where the implicit constant may depend on $\epsilon$.
\end{remark}}
While many of the questions about the M\"obius function are out of reach with current methods, researchers have attempted to study similar questions when $\mu$ is replaced by the so-called \textit{Random Multiplicative Functions}.

There are two main models of random multiplicative functions: 
\begin{itemize}
    \item \textbf{Rademacher random multiplicative functions:} let $(f(p))_{p\in\P}$ denote a sequence of i.i.d. random variables indexed by the primes and taking the values $\pm 1$ with equal probability. A Rademacher random multiplicative function is a sequence of random variables $(f(n))_{n\in \N}$ defined multiplicatively by\footnote{For convenience, we can assume that $f(n)=f(-n)$ for $n< 0$ and $f(0)=0$.} \[f(n)=\begin{cases}
        \prod\limits_{p\mid n} f(p)&\text{ if $n$ is square-free}\\
        0&\text{otherwise.}
    \end{cases}\]
    The idea is that this function might model the behaviour of the M\"obius function, which corresponds to the case when $f(p)=-1$ for all primes $p$.\footnote{
    We can also consider the so-called completely multiplicative Rademacher RMF defined by
    \[f(n)=
        \prod\limits_{p^\alpha\parallel n} f(p)^\alpha,\]
        which might model the Liouville function. However, we will focus our attention on the multiplicative case in this paper.}
    \item \textbf{Steinhaus random multiplicative functions} are defined in a similar way to Rademacher RMF, but with $f(p)$ distributed uniformly on the unit circle and with $f(n)$ being completely multiplicative. These random multiplicative functions are meant to model the Archimedean characters $n^{it}$, $t$ real.
\end{itemize}

In this paper, we will study random analogues of the Chowla conjecture, extending recent work of Klurman--Shkredov--Xu (\cite{Klurman2023}), Chinis--Shala (\cite{chinis2025random}) and de la Bret\`eche--Wang--Xu (\cite{delabreteche}). Also, we will study random analogues of the Goldbach conjecture motivated by the recent work of Mangerel (\cite{mangerel2024goldbach},\cite{mangerel2024shusterman}).
\subsection{Random Chowla Conjecture}

The Chowla conjecture (\cite{Chowla}) states that for the Liouville function $\lambda$ (or the M\"obius function $\mu$) and any polynomial $P(x)\in\Z[x]$, which is not of the form $P(x)=cg^2(x)$ for some
$g\in\Z[x]$, 
\[\sum_{n\leqslant x} \lambda(P(n))=o(x).\]

A natural random analogue of this conjecture arises when we replace $\lambda$ by a random multiplicative function. The case when $\lambda$ is replaced by a Steinhaus RMF was fully resolved by Klurman--Shkredov--Xu (\cite{Klurman2023}), who showed the following:
\begin{thm}
    Let $f$ be a Steinhaus random multiplicative function. Then for any polynomial $P(x)$ in $\mathbb{Z}[x]$ with $\deg P \geqslant 2$ which is not of the form $P(x) = w(x+c)^d$ for some $w \in \mathbb{Z}$, $c \in \mathbb{Q}$, as $N \to \infty$,\[
\frac{1}{\sqrt{N}} \sum_{n \leqslant N} f(P(n)) \xrightarrow{d} \mathcal{N}_\C(0,1),
\]
where $\xrightarrow{d}$ stands for convergence in distribution and $\mathcal{N}_\C$ stands for the standard complex Gaussian distribution. 
\end{thm}

The case when $\lambda$ is replaced by a Rademacher random multiplicative function was considered by Chinis--Shala (\cite{chinis2025random}), where they proved the following result:
\begin{thm}\label{thm:chinis-shala}
Let $f$ be a Rademacher RMF and let $P \in \mathbb{Z}[x]$ be a product of at least two distinct linear factors over $\mathbb{Z}$ or irreducible of degree $2$. Suppose that there is no prime $p$ for which $p^2 \mid P(n)$ for all $n \in \mathbb{N}$. Then, there exists a constant $\kappa_P > 0$ such that
\[
\frac{1}{\sqrt{\kappa_P N}} \sum_{n \leqslant N} f(P(n)) \xrightarrow{d} \mathcal{N}(0,1),
\]
as $N \to \infty$; that is, the partial sums $\sum_{n \leqslant N} f(P(n))$, when appropriately normalised, converge in distribution to a standard (real) Gaussian, as $N \to \infty$.
\end{thm}

\blue{Very recently, this was generalised by de la Bret\`eche--Wang--Xu (\cite{delabreteche})\footnote{\blue{This preprint appeared during the final stages of manuscript preparation; the research presented here was conducted entirely independently.}} to any polynomial having a positive proportion of square-free values. This class includes all admissible polynomials (meaning that the polynomial is square-free and there is no prime $p$ such that $p^2 \mid P(n)$ for all $n \in \mathbb{N}$) of degree up to $3$, and, conditional on the $abc$-conjecture, admissible polynomials of any degree.}

\begin{thm}[Corollary of Theorem 1.1 from \cite{delabreteche}]\label{thm:chowla-main}
Let $P$ be a non-linear square-free polynomial with integer coefficients and a positive leading coefficient. Moreover, suppose that there is a positive constant $\kappa_P$ such that
\[
\lim_{N\rightarrow\infty}
\frac{\#\{n\in\sqbra{N}:P(n)\text{ is square-free}\}}{N}
=\kappa_P.
\]

Let $f$ be a Rademacher RMF. Then
\[
\frac{1}{\sqrt{\kappa_PN}}\sum_{n\leqslant N}f(P(n))
\xrightarrow[]{d}\mathcal{N}(0,1).
\]
\end{thm}

We should also note that while Klurman--Shkredov--Xu, Chinis--Shala and de la Bret\`eche--Wang--Xu reduced the problem to computing the second and fourth moments of partial sums of RMFs at a polynomial argument and then applying McLeish’s martingale central limit theorem, \blue{we present an alternative proof of Theorem~\ref{thm:chowla-main} that is entirely independent of this martingale framework. }  Instead, we will directly estimate all moments of partial sums of Rademacher RMFs at polynomial arguments using a combinatorial argument. In particular, we will prove the following stronger statement, which is directly connected to moments of partial sums of Rademacher RMF at polynomial arguments:
\begin{thm}\label{thm:count-chowla-moments}
Let $n\in\N$ be fixed, and let $P$ be a square-free polynomial of degree $d$ at least $2$ with integer
coefficients \blue{and a positive leading coefficient}. Then there exists a constant $\mu=\mu(n,d)>0$, depending only on $d$
and $n$, such that the number of $n$-tuples $(x_1,\dots,x_n)$ such that
\begin{itemize}
    \item $x_i\in[N]$,
    \item $P(x_i)$ are square-free, pairwise distinct, and non-zero,
    \item $P(x_1)P(x_2)\cdots P(x_n)=\square$
\end{itemize}
is $O_{P,n,\ep}(N^{\frac{n}{2}-\mu+\ep})$.
\end{thm}
\blue{
\begin{remark}
    With a more delicate analysis, we can bound the number of solutions by $\exp\bra{Cn^2\cdot\bra{\frac{\log N}{\log\log N}+\log n}}N^{\frac{n}{2}-\frac{1}{40n}}$, where $C$ depends only on the polynomial $P$, and $N$ is sufficiently large in terms of $P$. Note that the first term in the exponent comes from the explicit maximal order of the divisor function, and the second comes from bounding the number of possible matrices $\tilde{A}$ (which will be defined below).
\end{remark}}
\begin{remark}\label{remark:count-chowla-moments}
The statement is trivial for $n=1$. Also, we can replace the condition \emph{$P(x_i)$ are distinct} by
\emph{$P(x_i)$ do not split into equal pairs} by taking the smallest $\mu$ among all
smaller values of $n$ for fixed $d$.
\end{remark}
\blue{This can be viewed as a Rademacher analogue of the paucity phenomenon established for moments of partial sums of Steinhaus RMFs by Wang--Xu (\cite{WangXu}). We devote Section~\ref{section:chowla} to this proof.}

\subsection{Random Goldbach Conjecture}

The Goldbach conjecture states that for any even number $n>2$, there exist prime numbers $p$ and $q$ such that $p+q=n$. This conjecture is also far beyond the reach of current methods. One of the main obstructions here is the so-called \textit{``parity problem''}. The famous result of Chen (\cite{Chen1973}) states that any sufficiently large even number can be represented as a sum of a prime and a product of at most two primes. 

On the other hand, the Goldbach conjecture trivially implies that for any even $N>2$, there exists an $n<N$ such that $\lambda(n)=\lambda(N-n)=-1$, and even this question is still open. Note that Chen's theorem does not tell us anything about $\lambda(N-n)$. A similar problem was recently addressed by Mangerel (\cite{mangerel2024goldbach}), who showed that for any sufficiently large $N$ there exists at least one $n<N$ such that $\lambda(n)=-1$ and $\lambda(N-n)=1$. Also, in a later paper (\cite{mangerel2024shusterman}), Mangerel showed that under the GRH, for any sufficiently large even integer $N$ there exists at least one $n<N$ such that $\lambda(n)=\lambda(N-n)=-1$.

In this paper, we will study similar problems when we replace $\lambda$ by a random multiplicative function. \blue{For instance, if $f$ is a Rademacher RMF, one might ask how often $f(n)$ and $f(N-n)$ take a specific pair of values. Alternatively, if $f$ is a Steinhaus RMF, which takes values continuously on the unit circle, we can study how often $f(n)$ and $f(N-n)$ lie in specific regions. For example, we might ask how often $f(n) \in I$ and $f(N-n) \in J$, where $I$ and $J$ are finite unions of arcs.}

For a Steinhaus random multiplicative function $f$, we will prove the following main result:

\begin{thm}\label{thm:main-steinhaus}
Let $F:\T\times \T\to\R$ be a measurable real-valued function of harmonic bounded variation (HBV)\footnote{\blue{See Definition~\ref{def:HBV} in Section~\ref{section:goldbach-steinhaus}.}}. Let $X$ and $Y$ be independent random variables, each distributed uniformly on $\T$. Define $G:\T\times \T\to\R$ by
\[
G(x,y)\coloneqq F(x,y)-\EE\bigl[F(x,Y)\bigr]-\EE\bigl[F(X,y)\bigr]+\EE\bigl[F(X,Y)\bigr].
\]
Set
\[
C\coloneqq \var\bigl(G(X,Y)\bigr)+\cov\bigl(G(X,Y),\,G(Y,X)\bigr),
\]
and assume that $C>0$. Let $f$ be a Steinhaus random multiplicative function. Then, as $N\to\infty$,
\[
\frac{\sum_{n<N} F\bigl(f(n),f(N-n)\bigr)-N\,\EE\bigl[F(X,Y)\bigr]}{\sqrt{CN}}
\ \xrightarrow{\,d\,}\ \mathcal{N}(0,1).
\]
\end{thm}
\begin{remark}
In this paper, we will mean by $\T$ both the unit circle in the complex plane and real numbers modulo $1$, depending on the context.
\end{remark}
\begin{remark}\label{rem:C-non-zero}
The constant $C$ is nonnegative. Moreover, the condition $C>0$ is equivalent to the assertion that the relation
\[
G(x,y)=-G(y,x)
\]
fails on a set of positive measure (with respect to Haar measure on $\T\times \T$). In the complementary case, as we shall see later, one obtains better than square-root cancellation.
\end{remark}

\blue{In particular, as a corollary, we can answer the question described above:}
\begin{corollary}\label{cor:product-indicators}
    Let $I$ and $J$ be unions of finitely many arcs on the unit circle. Let $a$ be the probability that a uniformly chosen point lies in $I$, let $b$ be the probability that it lies in $J$, and let $c$ be the probability that it lies in $I\cap J$. Assume moreover that $0<a,b<1$. Let $f$ be a Steinhaus random multiplicative function. Then, as $N\rightarrow \infty$
    \[
    \frac{\sum_{n<N}\mathbbm{1}(f(n)\in I)\mathbbm{1}(f(N-n)\in J)-abN}{\sqrt{\bra{a(1-a)b(1-b)+(c-ab)^2}N}}\xrightarrow{d}\mathcal{N}(0,1).
    \]
\end{corollary}
\blue{Another consequence of Theorem~\ref{thm:main-steinhaus} is the following result:}
\begin{corollary}\label{cor:complex-gaussian}
Let $f$ be a Steinhaus random multiplicative function. Then, as $N\rightarrow \infty$
    \[
    \frac{\sum_{n<N}f(n)f(N-n)}{\sqrt{2N}}\xrightarrow{d}\mathcal{N}_{\C}(0,1),
    \]
    where $\mathcal{N}_{\C}(0,1)$ denotes the standard complex Gaussian distribution with mean $0$ and variance $1$.
\end{corollary}

\blue{Finally, at the end of the paper, we will prove the following analogue of Corollary~\ref{cor:complex-gaussian} for the case where $f$ is a Rademacher RMF}:
\begin{thm}\label{thm:main-rademacher}
Let $f$ be a Rademacher random multiplicative function. Let $S(N)$ be the
number of positive integers $n<N$ such that $n$ and $N-n$ are both square-free.
Then, as $N\to\infty$,
\[
\frac{\sum_{n<N} f(n)f(N-n)}{\sqrt{2S(N)}}
\ \xrightarrow{\,d\,}\ \mathcal{N}(0,1).
\]
\end{thm}
\blue{From this, we deduce the following partial result, which serves as an analogue of Corollary~\ref{cor:product-indicators}}:
\begin{corollary}
\label{cor:rademacher-indicators}
Let $f$ be a Rademacher random multiplicative function. Let $S(N)$ be the
number of positive integers $n<N$ such that $n$ and $N-n$ are both square-free.
Suppose $c_1,c_2\in\{\pm1\}$ \emph{are such that} $c_1c_2=-1$. Then, as $N\to\infty$,
\[
\frac{\sum_{n<N}\mathbbm{1}(f(n)=c_1)\mathbbm{1}(f(N-n)=c_2)-\frac{S(N)}{4}}
{\sqrt{\frac{S(N)}8}}
\ \xrightarrow{\,d\,}\ \mathcal{N}(0,1).
\]
\end{corollary}\blue{
\begin{remark}
Using standard techniques (as in the proof of Theorem~3.1 in
\cite{lester_rudnick_2015}), we can show that
$S(N)=\rho(N)N+O_\ep(N^{\frac{2}{3}+\ep})$, where $\rho(N)=\prod_{p:p^2\nmid N}\Bigl(1-\frac{2}{p^2}\Bigr)\prod_{p:p^2\mid N}\Bigl(1-\frac{1}{p^2}\Bigr)$.
Note, in particular, that $\rho(N)\geqslant \prod_{p}\Bigl(1-\frac{2}{p^2}\Bigr)>0$,
and so $S(N)\asymp N$.

Therefore, we may replace $S(N)$ by $\rho(N)N$ in the denominator. The same applies to Theorem~\ref{thm:main-rademacher}. However, in the centring term $\frac{S(N)}{4}$ in the numerator, we cannot replace $S(N)$ by $\rho(N)N$, because the known error term is larger than $\sqrt{N}$. 
\end{remark}}
In the proofs of both Theorem~\ref{thm:main-steinhaus} and Theorem~\ref{thm:main-rademacher}, we will rely on McLeish’s martingale central limit theorem.
\subsection{Outline of the paper}
Our paper is organised as follows. In Section~\ref{section:chowla}, we prove Theorem~\ref{thm:chowla-main} by estimating all the moments of the partial sums of Rademacher RMF at polynomial arguments. In Section~\ref{section:goldbach-steinhaus}, we prove Theorem~\ref{thm:main-steinhaus} and deduce Corollaries \ref{cor:product-indicators} and \ref{cor:complex-gaussian}. This section, in fact, consists of two main parts: reducing the result to computing the second and fourth moments using McLeish's theorem, and the computation of the moments themselves, which requires delicate analysis. Finally, in Section~\ref{section:goldbach-rademacher}, we prove Theorem~\ref{thm:main-rademacher}, again utilising McLeish's theorem.
\section{Chowla's conjecture for Rademacher Random Multiplicative Functions}\label{section:chowla}

In this section, we will prove Theorem~\ref{thm:count-chowla-moments}, and show how it implies Theorem~\ref{thm:chowla-main}:
\begin{proof}[Proof of Theorem~\ref{thm:chowla-main}]
Since there are finitely many values of $n$ for which $P(n)\leqslant 0$, it is enough to consider $\frac{1}{\sqrt{\kappa_PN}}\sum_{\substack{n\leqslant N\\P(n)>0}}f(P(n))$. Let
\[
M_N^{(k)}
=
\EE\sqbra{
\bra{\frac{1}{\sqrt{\kappa_PN}}\sum_{\substack{n\leqslant N\\P(n)>0}}f(P(n))}^k
},
\]
the $k$-th moment of the $N$-th partial sum. It is enough to show that
\[
M^{(k)}_N\rightarrow
\begin{cases}
0,&k\text{ is odd},\\
(k-1)\doublefactorial,&k\text{ is even}.
\end{cases}
\tag{$\dagger$}
\]

Note that
\[
M_N^{(k)}
=\frac{\#\{(x_1,\dots,x_k)\in[N]^k:
P(x_i)>0\text{ is square-free and }
P(x_1)P(x_2)\cdots P(x_k)=\square \}}{\bra{\kappa_PN}^{\frac{k}{2}}}.
\]

By Theorem~\ref{thm:count-chowla-moments} and Remark~\ref{remark:count-chowla-moments}, the number of non-trivial $k$-tuples is
$o(N^{\frac{k}{2}})$, so they contribute $o(1)$ to $M_N^{(k)}$.

If $k$ is odd, there are no trivial solutions, and therefore $(\dagger)$ holds.

Now suppose $k$ is even. There are $(k-1)\doublefactorial$ ways to split the indices into pairs.
For a fixed pairing, the corresponding trivial contribution is
\[
\Bigl(
\#\{(x,y)\in[N]^2: P(x)=P(y)>0\text{ is square-free}\}
\Bigr)^{k/2}.
\]
Only $O_P(1)$ values of $P$ have more than one preimage or are non-positive, and therefore
\[
\#\{(x,y)\in[N]^2: P(x)=P(y)>0\text{ is square-free}\}
=\kappa_PN+o(N).
\]
Hence, each pairing contributes
\[
\frac{(\kappa_PN+o(N))^{k/2}}{(\kappa_PN)^{k/2}}=1+o(1),
\]
and so
\[
M_N^{(k)}=(k-1)\doublefactorial+o(1).
\]

Finally, note that in the previous count, we may have overcounted tuples for which two
distinct pairs take the same value, but the number of such tuples is again
negligible, so the conclusion remains unchanged.
\end{proof}
\subsection{Proof of Theorem~\ref{thm:count-chowla-moments}}
\blue{Generally speaking, both Wang--Xu and Chinis--Shala used the following approach for computing moments of partial sums of random multiplicative functions: manually considering multiple cases for the values of $x_i$ or $\gcd(P(x_i), P(x_j))$ for different $i$ and $j$, and finding an improvement in each case. It is worth noting that in all of these cases, the bounds between them are of the form $N^c$, where $c$ is a fixed constant that depends only on the degree of $P$. We will proceed with a similar approach. However, instead of manually considering different cases and selecting the boundaries, we will consider the cases more systematically. We describe the specific situations (Lemmas \ref{lem:improvement1}--\ref{lem:improvement4}) where we can get a power saving. Finally, we show that these four families cover all the cases. As in the work of Wang--Xu and Chinis--Shala, the main number-theoretic inputs in this mostly combinatorial argument are the divisor bound and the following two technical lemmata:}
\begin{lemma}\label{lem:bound-divisible}
Let $P$ be as above and $1\leqslant l\leqslant N^{O(1)}$. Then
\[
\#\{x\in\sqbra{N}:l\mid P(x)\}
\leqslant O_{P,\ep}\bra{\frac{N^{1+\ep}}{\min\{N,l\}}}.
\]
\end{lemma}
\begin{proof}
This is an almost immediate corollary of Proposition~2.2 from \cite{WangXu}, which states
that (note that $Q=P$ and $e_P=1$ in our setting)
\[
\#\{x\in\sqbra{N}:l\mid P(x)\}
\leqslant d^{\omega(l)}\abs{\Delta_P}^{\frac{1}{2}}\cdot \bra{1+\frac{N}{l}}.
\]

Now we only need to note that
\[
d^{\omega(l)}=O_\ep(l^\ep)=O_\ep(N^\ep),
\]
since $l\leqslant N^{O(1)}$.
\end{proof}
\begin{lemma}\label{lem:bound-big-gcd}
Let $P$ be as above, and let $a$, $b$ be two distinct non-zero integers. Then the
number of $(x,y)\in \sqbra{N}\times\sqbra{N}$ such that
\[
aP(x)-bP(y)=0
\]
is $O_\ep(N^{\frac{1}{2}+\ep})$.
\end{lemma}
\begin{proof}
Without loss of generality, assume that $a>0$. Note that if now $b\leqslant 0$, then $P(x)$ or $P(y)$ must be non-positive, and
there are only $O_P(1)$ positive integer points where $P$ is non-positive. Once one of $x,y$ is fixed, there are at most
$d$ possible values for the other, so there are $O_P(1)$ solutions in total in this case.

If $a,b>0$, then this follows directly from the proof of Lemma~2.7 in
\cite{WangXu}.
\end{proof}
\blue{
\begin{remark}
    We can replace $\frac{1}{2}$ by $\frac{1}{3}$ in this lemma, as in Lemma~2.5 of \cite{chinis2025random}. However, this is not necessary for our proof.
\end{remark}
}

\blue{We start with the following lemma, which helps us keep track of the greatest common divisors of $P(x_i)$ and $P(x_j)$ for different $i$ and $j$.}
\begin{lem}
Let $n\geqslant 2$ and $(a_1,\dots,a_n)$ be an $n$-tuple of square-free positive
integers whose product is a perfect square. Then there exists an $n\times n$ matrix
$A$ with positive integer entries, and with the following properties:
\begin{itemize}
    \item $A$ is symmetric,
    \item $A_{ii}=1$ for all $i$,
    \item for all $i$,
    $\prod_j A_{ij}=\prod_j A_{ji}=a_i$.
\end{itemize}
\end{lem}
\begin{remark}
We will call $A$ a $\gcd$-table, and we will apply the lemma to the tuple
$(\abs{P(x_1)},\dots,\abs{P(x_n)})$. The motivation is that $A_{ij}$ should be the
$\gcd(P(x_i),P(x_j))$, however, this is not quite enough to ensure that the product
in each row and column is exactly a value of $P$, and we have to \emph{reduce} the entries of $A$.

Moreover, note that the existence of such a matrix encodes the fact that the product
is a perfect square.
\end{remark}
\begin{proof}
We argue by induction on $n$. For $n=2$, the statement is trivial, since $a_1=a_2$,
and so we can take
\[
A=
\begin{pmatrix}
1 & a_1 \\
a_1 & 1
\end{pmatrix}.
\]

Now set $A_{1n}=A_{n1}=\gcd(a_1,a_n)$ and then define
\[
A_{in}=A_{ni}=\gcd\bra{a_i,\frac{a_n}{\prod_{j<i}{A_{jn}}}}
\]
for $i<n$, and $A_{nn}=1$.

Note that for all $i$, $A_{ni}=A_{in}\mid a_i$. Moreover,
\[
\prod_{i<n}A_{in}
=\prod_{i<n}\gcd\bra{a_i,\frac{a_n}{\prod_{j<i}{A_{jn}}}}
=\gcd(a_n,a_1\cdots a_{n-1})
=a_n,
\]
where the last equality comes from the fact that $a_n$ is square-free and
$a_1\cdots a_{n-1}a_n$ is a perfect square.

Now apply the induction hypothesis to $(a_1',\dots,a_{n-1}')$, where
$a_i'=\frac{a_i}{A_{in}}$. Clearly, all $a_i'$ are square-free integers, and
moreover
\[
a_1'\cdots a_{n-1}'
=\frac{a_1\cdots a_n}{a_n\prod_i A_{in}}
=\frac{a_1\cdots a_n}{a_n^2},
\]
so $a_1'\cdots a_{n-1}'$ is a perfect square. Therefore, there exists a matrix $A'$
satisfying all the required conditions for $(a_1',\dots,a_{n-1}')$, and gluing to it
the entries $A_{in},A_{ni}$ defined above produces a valid matrix for
$(a_1,\dots,a_n)$. This completes the induction.
\end{proof}

Note that for each matrix as in the lemma, we can have at most $O_{P,n}(1)$ tuples
$(x_1,\dots,x_n)$ such that it is a valid $\gcd$-matrix for
$(\abs{P(x_1)},\dots,\abs{P(x_n)})$. Hence, it is enough to bound the number of valid
$\gcd$-matrices, where the product in each row is the absolute value of $P$,
evaluated at a positive integer up to $N$.

Let now $M$ be a large constant, depending on $n$ and $d$, for example $
M=(10^6nd)^{100}
$ will work.
Consider a symmetric matrix $\tilde{A}$ defined by $
\tilde{A}_{ij}=M^{-1}\cdot\lfloor M\log_NA_{ij}\rfloor$.
Note that $\prod_j A_{ij}\ll_P N^d$, and therefore $\sum_j\tilde{A}_{ij}\leqslant d$ for all sufficiently large $N$.

Then $\tilde{A}$ has at most $(dM+1)^{n^2}=O_{n,d}(1)$ possible values, so it is enough to bound the number of solutions for each fixed $\tilde{A}$\blue{, i.e. we will distinguish cases according to $\tilde{A}$}.

In general, we reconstruct the values of $x_1,\dots,x_n$ by the following algorithm
(this procedure will often fail to produce a valid solution, but every valid solution
can be obtained in this way):
\begin{itemize}
    \item Choose $x_1$.
    \item Choose $A_{1i}$ for all $i$ so that their product is $\abs{P(x_1)}$ and
    they lie in the intervals defined by $\tilde{A}$.
    \item Choose $x_2$ such that $P(x_2)$ is divisible by $A_{12}$.
    \item Choose $A_{2i}$ as above.
    \item Choose $x_3$ such that $P(x_3)$ is divisible by $A_{13}A_{23}$.
    \item Choose $A_{3i}$.
    \item Continue in this way until all the $x_i$ are chosen.
\end{itemize}
    Moreover, we may permute the rows and columns beforehand.  

    Define
$
\beta_i=\frac{1}{d}\sum_j\tilde{A}_{ij}\leqslant 1
$.
Then, \[\abs{x_i}^d\ll_P\abs{P(x_i)}=\prod_j A_{ij}\leqslant N^{\sum_j\bra{\tilde{A}_{ij}+\frac{1}{M}}}\leqslant N^{d\beta_i+\frac{n}{M}},\]
so
$\abs{x_i}\ll_P N^{\beta_i+\frac{n}{M}}$.
By Lemma~\ref{lem:bound-divisible}, as $A_{1i}\cdots A_{i-1,i}
<
N^{\sum_{j<i}\tilde{A}_{ji}+\frac{n}{M}}\leqslant N^{O_{n,d}(1)},$ we have
$
O_\ep\bra{
\frac{N^{\beta_i+\frac{n}{M}+\ep}}
{\min\{N^{\beta_i+\frac{n}{M}},A_{1i}\cdots A_{i-1,i}\}}
}
$ possibilities for $x_i$.
    
Also, $N^{\sum_{j<i}\tilde{A}_{ji}}
\leqslant
A_{1i}\cdots A_{i-1,i}$. Therefore, we obtain
\[
O_\ep\bra{
\frac{N^{\beta_i+\frac{n}{M}+\ep}}
{\min\{N^{\beta_i+\frac{n}{M}},N^{\sum_{j<i}\tilde{A}_{ji}}\}}
}
\]
possibilities for $x_i$.

Let $\gamma_i=\sum_{j<i}\tilde{A}_{ji}$.
Then, when choosing $x_i$, we have
$
O_\ep(N^{\beta_i-\gamma_i+\frac{n}{M}+\ep})
$ possibilities if $\gamma_i<\beta_i$, and
$
O_\ep(N^{\frac{n}{M}+\ep})
$ otherwise. Note that this introduces an ``error'' of size
$O_\ep(N^{\frac{n}{M}+\ep})$ at each step, which overall contributes
$O_\ep(N^{\frac{n^2}{M}+\ep})$.

Also, once $x_i$ has been chosen, there are $O_\ep(N^\ep)$ possibilities for choosing the
corresponding $A_{ij}$ by the divisor bound.

Let $I$ be the set of $i$ such that $\gamma_i<\beta_i$, and let $J$ be its
complement. Thus, the problem reduces to bounding
\[
\sum_{i\in I}(\beta_i-\gamma_i).
\]
    
    We now prove the following proposition, which plays a crucial role in our proof.
    \begin{prop}\label{prop:main-tech}
For all $t$ from $0$ to $n$ we have
\[
\sum_{I\cap\{i:i>t\}}(\beta_i-\gamma_i)
\leqslant
\frac{1}{2d}\sum_{i>t,j}\tilde{A}_{ij}
-\frac{d-1}{d}\sum_{I\cap\{i:i>t\}}\gamma_i
\leqslant
\frac{1}{2d}\sum_{i>t,j}\tilde{A}_{ij}
\leqslant
\frac{n-t}{2}.
\]
\end{prop}
\begin{proof}
First note that the last two inequalities are immediate, since for all $i$ we have
$\sum_j \tilde{A}_{ij}\leqslant d$ and $\gamma_i\geqslant0$.

Let
\[
I'=I\cap\{i:i>t\},
\qquad
J'=J\cap\{i:i>t\}.
\]
We consider two cases.

\begin{enumerate}
    \item If $\sum_{I'} \beta_i<S$, where $S$ will be specified later, then
    \[
    \sum_{i\in I'}(\beta_i-\gamma_i)
    \leqslant
    S-\sum_{i\in I'} \gamma_i.
    \]

    \item If $\sum_{I'}\beta_i\geqslant S$, then note that for all $j$,
    $\gamma_j\leqslant d\beta_j$, so
    $
    \sum_{J'}\gamma_j\leqslant d\sum_{J'}\beta_j.
    $

    Observe that
    \[
    2\sum_{i>t} \gamma_i
    =
    \sum_{i>t\text{ or }j>t}\tilde{A}_{ij}
    \geqslant
    \sum_{i>t,j}\tilde{A}_{ij},
    \]
    and hence
    \[
    \sum_{I'}\gamma_i
    \geqslant
    \frac{1}{2}\sum_{i>t,j}\tilde{A}_{ij}
    -d\sum_{J'}\beta_j.
    \]

    Therefore
    \begin{align*}
    \sum_{I'}(\beta_i-\gamma_i)
    &\leqslant
    \sum_{I'} \beta_i
    -\frac{1}{2}\sum_{i>t,j}\tilde{A}_{ij}
    +d\sum_{J'}\beta_j \\
    &=
    \frac{1}{2}\sum_{i>t,j}\tilde{A}_{ij}
    -(d-1)\sum_{I'}\beta_i \\
    &\leqslant
    \frac{1}{2}\sum_{i>t,j}\tilde{A}_{ij}
    -(d-1)S.
    \end{align*}
    Here we used that
    $d\sum_{i>t}\beta_i=\sum_{i>t,j}\tilde{A}_{ij}$.
\end{enumerate}

Finally, taking
$
S=\frac{1}{2d}\sum_{i>t,j}\tilde{A}_{ij}
+\frac{1}{d}\sum_{i\in I'}\gamma_i,
$
we obtain the desired bound in both cases.
\end{proof}

        Combining the bound from Proposition~\ref{prop:main-tech} with $t=0$ with the
argument above, we obtain the estimate
$
O_{P,\ep,n,d}(N^{\frac{n}{2}+\frac{n^2}{M}+\ep})
$
for the number of solutions corresponding to a fixed $\tilde{A}$.

This is already close to the desired bound, and we require only a slight improvement to achieve the final result. To achieve this, we prove several lemmas, each yielding an improvement in a particular case, and finally show that any case can be reduced to one of these lemmas.

\begin{lemma}\label{lem:improvement1}
Suppose there is some $i$ such that $\beta_i\leqslant 1-\frac{1}{10n}$. Then the
total number of solutions for such $\tilde{A}$ is bounded by
\[
O_{P,\ep,n,d}\bra{N^{\frac{n}{2}-\frac{1}{20n}+\frac{n^2}{M}+\ep}}.
\]
\end{lemma}
\begin{proof}
By Proposition~\ref{prop:main-tech} with $t=0$,
\[
\sum_I(\beta_i-\gamma_i)\leqslant\frac{1}{2d}\sum_{ij}\tilde{A}_{ij}.
\]
But in this case, we can estimate the right-hand side more sharply:
\[
\sum_I(\beta_i-\gamma_i)
\leqslant
\frac{1}{2d}\sum_{ij}\tilde{A}_{ij}
=
\frac{1}{2}\sum_i\beta_i
\leqslant
\frac{1}{2}\cdot \bra{\Bigl(1-\frac{1}{10n}\Bigr)+(n-1)}
=
\frac{n}{2}-\frac{1}{20n},
\]
which proves the lemma.
\end{proof}

From now on, we may assume that all $\beta_i$ are at least $1-\frac{1}{10n}$.
\begin{lemma}\label{lem:improvement2}
Suppose $\tilde{A}$ has an entry $\tilde{A}_{ij}$ which is less than $\beta_i$ and
greater than or equal to $\frac{1}{10n}$. Then the total number of solutions for such
$\tilde{A}$ is bounded by
\[
O_{P,\ep,n,d}\bra{N^{\frac{n}{2}-\frac{1}{20n}+\frac{n^2}{M}+\ep}}.
\]
\end{lemma}
\begin{proof}
Since simultaneously permuting the rows and columns of $A$ and $\tilde{A}$ does not
change the number of solutions, i.e. permuting the $(x_k)_{k=1}^n$, we may assume without loss of generality that this
entry is $\tilde{A}_{12}=\tilde{A}_{21}$, and so
\[
\frac{1}{10n}\leqslant \tilde{A}_{12} <\beta_2.
\]

Then, by definition, $2\in I$, and therefore
\[
\sum_{i\in I}\gamma_i\geqslant \gamma_2=\tilde{A}_{12}\geqslant \frac{1}{10n}.
\]
Applying the first inequality from Proposition~\ref{prop:main-tech}, we obtain
\[
\sum_{I}(\beta_i-\gamma_i)
\leqslant
\frac{1}{2d}\sum_{i,j}\tilde{A}_{ij}
-\frac{d-1}{d}\sum_{I}\gamma_i
\leqslant
\frac{n}{2}-\frac{d-1}{d\cdot(10n)}
\leqslant
\frac{n}{2}-\frac{1}{20n},
\]
which proves the lemma.
\end{proof}

Thus, from now on, we may further assume that every entry of $\tilde{A}$ is either
less than $\frac{1}{10n}$, or greater than or equal to $\beta_i$ for the corresponding
$i$. We will call such entries \emph{big}. By the previous assumption, all big
entries are at least $1-\frac{1}{10n}$. Moreover, by symmetry, if $\tilde{A}_{ij}$ is
big, then it must be at least $\beta_j$ as well.

    \begin{lemma}\label{lem:improvement3}
Suppose there exist $i\neq j$ such that both
\[
\sum_{k\neq j}\tilde{A}_{ik}\leqslant\frac{1}{5}
\qquad\text{and}\qquad
\sum_{k\neq i}\tilde{A}_{kj}\leqslant\frac{1}{5}.
\]
Then the total number of solutions for such $\tilde{A}$ is
bounded by
\[
O_{P,\ep,n,d}\bra{N^{\frac{n}{2}-\frac{1}{10}+\frac{2n}{M}+\frac{n^2}{M}+\ep}}.
\]
\end{lemma}
\begin{proof}
Without loss of generality, assume $i=1$ and $j=2$.

Instead of applying the algorithm above directly, we bound the number of possible
values of $x_1$ and $x_2$ simultaneously using Lemma~\ref{lem:bound-big-gcd}.

Indeed,
\[
\abs{\frac{P(x_1)}{P(x_2)}}
=\frac{\prod_{j>2} A_{1j}}{\prod_{i>2} A_{i2}}
=\frac{a}{b}.
\]

Note that
\[
a\leqslant N^{\frac{n}{M}+\sum_{j>2} \tilde{A}_{1j}},
\qquad
b\leqslant N^{\frac{n}{M}+\sum_{i>2} \tilde{A}_{i2}}.
\]
Therefore, by the assumptions of the lemma,
$
a,b\leqslant N^{\frac{1}{5}+\frac{n}{M}},
$
and so there are $O\bra{N^{\frac{2}{5}+\frac{2n}{M}}}$ possibilities for $\frac{P(x_1)}{P(x_2)}$.

Moreover, we assume that $P(x_1)\neq P(x_2)$ in Theorem~\ref{thm:count-chowla-moments}. Therefore,
Lemma~\ref{lem:bound-big-gcd} implies that for each possible value of
$\frac{P(x_1)}{P(x_2)}$ there are $
\ll_\ep N^{\frac{1}{2}+\ep}
$
possible solutions. Hence, in total, there are $
\ll_\ep N^{1-\frac{1}{10}+\frac{2n}{M}+\ep}$
possibilities for the pair $(x_1,x_2)$.

After choosing $x_1$ and $x_2$, we continue with the usual algorithm. By
Proposition~\ref{prop:main-tech} for $t=2$,
\[
\sum_{i\in I;i>2}(\beta_i-\gamma_i)\leqslant\frac{n-2}{2}.
\]
Therefore, there are $
O_\ep\bra{N^{\frac{n-2}{2}+\frac{n^2}{M}+\ep}}
$
ways to complete the table after $x_1$ and $x_2$ are fixed. Combining this with the
previous estimate proves the lemma.
\end{proof}

    Finally, we have the following.
    \begin{lemma}\label{lem:improvement4}
Suppose there is a row in $\tilde{A}$ which contains at least two \emph{big}
entries. Then the total number of solutions for such $\tilde{A}$ is bounded by
\[
O_{P,\ep,n,d}\bra{N^{\frac{n}{2}-\frac{1}{2}+\frac{n^2}{M}+\ep}}.
\]
\end{lemma}
\begin{proof}
Without loss of generality, we may assume this is the first row, and the two big
entries are $\tilde{A}_{12}$ and $\tilde{A}_{13}$.

Then clearly $2$ and $3$ are not in $I$ by definition. Therefore, applying
Proposition~\ref{prop:main-tech} with $t=3$, we obtain
\[
\sum_I(\beta_i-\gamma_i)
\leqslant
\beta_1+\sum_{i\in I;i>3}(\beta_i-\gamma_i)
\leqslant
\beta_1+\frac{n-3}{2}
\leqslant
\frac{n-1}{2},
\]
which gives the stated bound.
\end{proof}

    We can now combine the previous lemmas to complete the proof. By
Lemma~\ref{lem:improvement1}, we know that
$
\beta_1\geqslant 1-\frac{1}{10n}\geqslant\frac{1}{2}$.
Hence, $
\sum_i\tilde{A}_{1i}\geqslant \frac{d}{2}$,
so for some $i$ we must have $
\tilde{A}_{1i}\geqslant \frac{1}{2n}$.
Without loss of generality, this is $\tilde{A}_{12}$.
    
By Lemma~\ref{lem:improvement2}, we may assume
$
\tilde{A}_{12}\geqslant \beta_1,\beta_2.
$
Now apply Lemma~\ref{lem:improvement3} with $i=1$ and $j=2$. Then either
\[
\sum_{j\neq 2}\tilde{A}_{1j}\geqslant \frac{1}{5}
\qquad\text{or}\qquad
\sum_{i\neq 1}\tilde{A}_{i2}\geqslant \frac{1}{5}.
\]
By symmetry, we may assume $
\sum_{j\neq 2}\tilde{A}_{1j}\geqslant \frac{1}{5}$.
Hence, for some $j\neq 2$, we have
$
\tilde{A}_{1j}\geqslant \frac{1}{5(n-1)}>\frac{1}{10n}$, and so by Lemma~\ref{lem:improvement2},
$
\tilde{A}_{1j}\geqslant \beta_1,\beta_j$.
    
    Finally, we apply Lemma~\ref{lem:improvement4} to the first row, using the two big
entries $\tilde{A}_{12}$ and $\tilde{A}_{1j}$. Therefore, every possible $\tilde{A}$
is covered by one of the above lemmas, and hence the total number of solutions for a
fixed $\tilde{A}$ is bounded by $
O_{P,\ep,n}\bra{N^{\frac{n}{2}-\frac{1}{20n}+\frac{3n^2}{M}+\ep}}$.
This bound is weaker than the one in each individual improvement lemma, and is
therefore valid uniformly.

    By the choice of $M$ above, $\frac{3n^2}{M}<\frac{1}{40n}$, and we may take $\mu=\frac{1}{40n}$ in the theorem. This completes the proof.
\section{Goldbach's conjecture for Steinhaus Random Multiplicative Functions}\label{section:goldbach-steinhaus}
\subsection{Overview of the main result}
We recall the standard definition of two-dimensional harmonic bounded variation for periodic functions, as used in \cite{saakyan} and \cite{waterman}.
\begin{defn}\label{def:HBV}
Let $I=(a,b)$ and $J=(c,d)$ be intervals.
Define increments
\[
  F(I,y)=F(b,y)-F(a,y),\quad F(x,J)=F(x,d)-F(x,c),
\]
\[
  F(I,J)=F(a,c)-F(a,d)-F(b,c)+F(b,d).
\]
Let $\Omega$ denote the family of finite collections of pairwise disjoint intervals in $[0,1)$.
For collections $E_1=\{I_i\}\in\Omega$ and $E_2=\{J_j\}\in\Omega$ write the indices so that $i=1,2,\dots$.
Define
\[
  HV_1(F)=\sup_{y\in[0,1)}\sup_{E_1\in\Omega}\sum_i \frac{|F(I_i,y)|}{i},
\]
\[
  HV_2(F)=\sup_{x\in[0,1)}\sup_{E_2\in\Omega}\sum_j \frac{|F(x,J_j)|}{j},
\]
\[
  HV_{1,2}(F)=\sup_{E_1,E_2\in\Omega}\sum_i\sum_j \frac{|F(I_i,J_j)|}{ij}.
\]
We say that a measurable double periodic function $F$ has \textbf{harmonic bounded variation}, denoted by $F\in\HBV(\T^2)$ if \[HV_1(F)+HV_2(F)+HV_{1,2}(F)<\infty.\]
\end{defn}
\begin{remark}
    Similarly, for a periodic function $F$ we write $F\in\HBV(\T)$ if \[HV(F)=\sup_{E\in\Omega}\sum_i \frac{|F(I_i)|}{i}<\infty.\]
\end{remark}
For a function $F \in HBV(\T^2)$, we denote the truncated rectangular Fourier partial sums by:
\[
  S_{N_1,N_2}(F;x,y) = \sum_{|n|\leqslant N_1} \sum_{|m|\leqslant N_2} c_{n,m} e^{2\pi i (nx+my)},
\]
where $c_{n,m}$ are the Fourier coefficients of $F$.
We will need the following result about HBV functions first established by Sahakian in \cite{saakyan}, which is explicitly mentioned as part of Theorem C in \cite{waterman}.
\begin{thm}\label{thm:HBV-bounded}
    If $F\in\HBV(\T^2)$, then the rectangular partial sums $S_{\boldsymbol{\cdot},\boldsymbol{\cdot}}(F;\boldsymbol{\cdot},\boldsymbol{\cdot})$ are uniformly bounded.
\end{thm}

We will also need the following elementary lemma:
\begin{lem}\label{lem:HBV-product}
    Let $h_1,h_2:\T\rightarrow \R$ be bounded, measurable, real-valued functions with bounded harmonic variation. Then the function $F$ defined by $F(x,y)=h_1(x)h_2(y)$ is in $\HBV(\T^2)$.
\end{lem}
\begin{proof}
    Note that $HV_1(F)\leqslant HV(h_1)\cdot \aabs{h_2}_\infty<\infty$. Similarly, $HV_2(F)<\infty$. Finally, observe that $F(I,J)=h_1(I)h_2(J)$. Hence, $HV_{1,2}(F)\leqslant HV(h_1)\cdot HV(h_2)<\infty$ which concludes the proof.
\end{proof}
\blue{We now demonstrate how Corollaries~\ref{cor:product-indicators} and \ref{cor:complex-gaussian} follow from our main result, Theorem~\ref{thm:main-steinhaus}.}
\begin{proof}[Proof of Corollary~\ref{cor:product-indicators}]
    Define $F=\mathbbm{1}(x\in I)\mathbbm{1}(y\in J)$. Since $I$ and $J$ are unions of finitely many arcs, the functions $\mathbbm{1}(x\in I)$ and $\mathbbm{1}(x\in J)$ are bounded and measurable with bounded variation (and hence bounded harmonic variation). Therefore, by Lemma~\ref{lem:HBV-product}, the function $F$ has HBV. The remaining hypotheses of Theorem~\ref{thm:main-steinhaus} are immediate.

It remains to compute $\EE[F(X,Y)]$ and the constant $C$ from Theorem~\ref{thm:main-steinhaus}:
\begin{align*}
\EE(F(X,Y)) &= \EE\mathbbm{1}(X\in I)\EE\mathbbm{1}(Y\in J)=ab.\\
\var(G(X,Y)) &= \EE\bigl(\mathbbm{1}(X\in I)-a\bigr)^2\bigl(\mathbbm{1}(Y\in J)-b\bigr)^2=a(1-a)b(1-b).\\
\cov(G(X,Y),G(Y,X)) &= \bigl(\EE(\mathbbm{1}(X\in I)-a)(\mathbbm{1}(X\in J)-b)\bigr)^2\\
&= (\EE\mathbbm{1}(X\in I\cap J)-ab)^2=(c-ab)^2.
\end{align*}

    Substituting these expressions into Theorem~\ref{thm:main-steinhaus} yields the stated conclusion. 
    
    Note that since $0<a,b<1$, we have $a(1-a)b(1-b)>0$, and therefore $C=a(1-a)b(1-b)+(c-ab)^2>0$.
\end{proof}
\begin{proof}[Proof of Corollary~\ref{cor:complex-gaussian}]
    We cannot apply Theorem~\ref{thm:main-steinhaus} directly with $F(x,y)=xy$, since $xy$ is not real-valued. Following the approach used in the preprint version of \cite{Klurman2023} (in the proof of Theorem~1.1), we instead consider the real-valued functions
    \[
    F_1=\Re(xy)\qquad\text{and}\qquad F_2=\Im(xy).
    \]
    It suffices to show that, for any real $a$ and $b$ with $(a,b)\neq(0,0)$,
    \[
    \frac{\sum_{n<N}\bigl(aF_1(f(n),f(N-n))+bF_2(f(n),f(N-n))\bigr)}{\sqrt{\frac{(a^2+b^2)}{2}\cdot 2N}}
    \xrightarrow{d}\mathcal{N}(0,1).
    \]
    Indeed, this implies that the limiting distribution of
    $\sum_{n<N}f(n)f(N-n)/\sqrt{2N}$ is $\mathcal{N}_{\mathbb{C}}(0,1)$.

    Note that
    \[
    F_1(x,y)=\Re x\,\Re y-\Im x\,\Im y,
    \]
    and each term on the right-hand side has HBV by Lemma~\ref{lem:HBV-product}; hence $F_1$ has HBV. The same argument applies to $F_2$. Therefore, Theorem~\ref{thm:main-steinhaus} applies to $F\coloneqq aF_1+bF_2$ provided the corresponding constant $C$ is non-zero; in fact, $C=a^2+b^2\neq 0$.

    Let $G$ be defined from $F$ as in Theorem~\ref{thm:main-steinhaus}. Then:
    \begin{itemize}
        \item $\EE\bra{F(X,Y)}=0$.
        \item $\var(G(X,Y))=\EE\bra{aF_1+bF_2}^2=\EE\bra{a\Re(XY)+b\Im(XY)}^2$.
        \item $\cov(G(X,Y),G(Y,X))=\var(G(X,Y))$ since $F$ is symmetric.
    \end{itemize}
    Writing $Z=XY$, we note that $Z$ is uniform on the unit circle, i.e.\ $Z=e^{i\Theta}$ with $\Theta\sim\mathrm{U}[0;2\pi]$. Hence
    \[
    a\Re(XY)+b\Im(XY)=a\cos\Theta+b\sin\Theta=\sqrt{a^2+b^2}\cos(\Theta+\theta_0),
    \]
    for some $\theta_0$, and therefore
    \[
    \var(G(X,Y))=(a^2+b^2)\EE\cos^2(\Theta+\theta_0)=\frac{a^2+b^2}{2}.
    \]
    Consequently, Theorem~\ref{thm:main-steinhaus} yields the desired convergence, and the corollary follows.
\end{proof}

\blue{We are now ready to prove Theorem~\ref{thm:main-steinhaus}. As mentioned in the Introduction, the proof consists of two main parts: using McLeish's theorem to reduce the result to the computation of second and fourth moments, and then establishing those moment estimates directly.}
\subsection{Reduction of Theorem~\ref{thm:main-steinhaus} to moment estimates}
\subsubsection{Independence lemma}

\begin{defn}[Multiplicative independence and span]
We say that positive integers $n_1,\dots,n_k$ are \emph{(multiplicatively) independent}
if whenever $\beta_1,\dots,\beta_k\in\Q$ satisfy
\[
n_1^{\beta_1}\cdots n_k^{\beta_k}=1,
\]
then $\beta_1=\cdots=\beta_k=0$.

We also define the \emph{(multiplicative) span} of $n_1,\dots,n_k$ as
\[
\Span(n_1,\dots,n_k)
:=\left\{\prod_{i=1}^k n_i^{\beta_i}:\ \beta_i\in\Q\right\}\subset \R_{>0}.
\]
\end{defn}

\begin{remark}\label{rem:valuation-vectors}
For a prime $p$, write $v_p(n)$ for the $p$-adic valuation of $n$, and set
$v(n):=(v_p(n))_{p\ \mathrm{prime}}\in\bigoplus_p \Z$ (finite support).
Then $n_1,\dots,n_k$ are independent if and only if the vectors
$v(n_1),\dots,v(n_k)$ are linearly independent over $\Q$.
Equivalently, $n\in \Span(n_1,\dots,n_k)$ if and only if
$v(n)$ lies in the $\Q$-linear span of $\{v(n_i)\}_{i=1}^k$.
\end{remark}

\begin{lemma}\label{lem:independence-lemma-strong}
Let $n,n_1,\dots,n_k$ be positive integers. If $n\notin \Span(n_1,\dots,n_k)$, then
$f(n)$ is independent of the random vector $(f(n_1),\dots,f(n_k))$, where $f$ is a Steinhaus RMF.
\end{lemma}

\begin{proof}
Since $f(n)$ and $f(n_i)$ take values in $\T$, their joint law on
$\T^{k+1}$ is determined by its Fourier coefficients. Thus, it suffices to show that
for every choice of integers $m,m_1,\dots,m_k\in\Z$,
\begin{equation}\label{eq:fourier-factorization}
\EE\!\left[f(n)^m\prod_{i=1}^k f(n_i)^{m_i}\right]
=\EE[f(n)^m]\cdot \EE\!\left[\prod_{i=1}^k f(n_i)^{m_i}\right].
\end{equation}

Let $P$ be the finite set of primes dividing $n n_1\cdots n_k$, and enumerate it as
$P=\{p_1,\dots,p_L\}$. Write
\[
n=\prod_{j=1}^L p_j^{a_j},
\qquad
n_i=\prod_{j=1}^L p_j^{a_{i,j}}
\quad (1\le i\le k),
\]
where $a_j=v_{p_j}(n)$ and $a_{i,j}=v_{p_j}(n_i)$. By complete multiplicativity,
\[
f(n)^m\prod_{i=1}^k f(n_i)^{m_i}
=\prod_{j=1}^L f(p_j)^{\,m a_j+\sum_{i=1}^k m_i a_{i,j}}.
\]

Note that for every prime $p$ and every integer $r\in\Z$,
\begin{equation}\label{eq:steinhaus-moments}
\EE[f(p)^r]=
\begin{cases}
1, & r=0,\\
0, & r\neq 0.
\end{cases}
\end{equation}
Therefore, using independence over primes and \eqref{eq:steinhaus-moments}, we obtain
\begin{align*}
\EE\!\left[f(n)^m\prod_{i=1}^k f(n_i)^{m_i}\right]
&=\prod_{j=1}^L \EE\!\left[f(p_j)^{\,m a_j+\sum_{i=1}^k m_i a_{i,j}}\right]
=\\&=
\begin{cases}
1, & \text{if } m a_j+\sum_{i=1}^k m_i a_{i,j}=0\ \text{for all }j,\\
0, & \text{otherwise.}
\end{cases}
\end{align*}
In vector form, this condition is
\[
m\,v(n)+\sum_{i=1}^k m_i\,v(n_i)=0.
\]
If $m\neq 0$ and this holds, then
\[
v(n)=-\frac{1}{m}\sum_{i=1}^k m_i\,v(n_i)\in \Span\{v(n_1),\dots,v(n_k)\},
\]
which implies $n\in \Span(n_1,\dots,n_k)$ by Remark~\ref{rem:valuation-vectors},
contradicting the assumption. Hence, any solution must have $m=0$, and therefore
\[
\EE\!\left[f(n)^m\prod_{i=1}^k f(n_i)^{m_i}\right]
=
\begin{cases}
1, & \text{if } m=0 \text{ and } \sum_{i=1}^k m_i\,v(n_i)=0,\\
0, & \text{otherwise.}
\end{cases}
\]

On the other hand, since $n\neq 1$ (because $1\in \Span(n_1,\dots,n_k)$),
\[
\EE[f(n)^m]=
\begin{cases}
1, & m=0,\\
0, & m\neq 0,
\end{cases}
\qquad
\EE\!\left[\prod_{i=1}^k f(n_i)^{m_i}\right]
=
\begin{cases}
1, & \sum_{i=1}^k m_i\,v(n_i)=0,\\
0, & \text{otherwise.}
\end{cases}
\]
Multiplying these two expressions gives exactly the same case distinction as above,
so \eqref{eq:fourier-factorization} holds for all $m,m_1,\dots,m_k\in\Z$. This
proves the claimed independence.
\end{proof}
As an immediate corollary, we also obtain the following result:
\begin{lemma}[Independence lemma]\label{lem:independence-lemma}
If $n_1,\dots,n_k$ are independent, then $f(n_1),\dots,f(n_k)$ are i.i.d.\ and
uniform on $\T$.
\end{lemma}

\begin{defn}[Strong dependence]
We say that positive integers $n_1,\dots,n_k$ are \emph{strongly dependent} if they are
dependent, and moreover for every $j\in\{1,\dots,k\}$,
\[
n_j\in \Span(n_1,\dots,n_{j-1},n_{j+1},\dots,n_k).
\]
\end{defn}
\subsubsection{Truncation of the Fourier series}\label{subsection:truncation}

Let $(c_{n,m})$ denote the Fourier coefficients of $F$. Throughout this and the next subsection, we will assume that $G=F$, and in Section~\ref{subsection:deducing-steinhaus-main} we will deduce Theorem~\ref{thm:main-steinhaus} in the most general case. This condition implies that $c_{n,0}=c_{0,m}=0$, and for $z\in \T$
\[
\EE F(z,Y)=\EE F(X,z)=0.
\]
\begin{remark}
    As in Theorem~\ref{thm:main-steinhaus}, $X$ and $Y$ are two independent uniform random variables on $\T$.
\end{remark}

We will truncate the Fourier expansion at
\[
t:=\lfloor\log\log\log N\rfloor,
\]
i.e.\ we will consider separately the Fourier coefficients with indices up to $t$ and all the others.

\begin{remark}
We only require $t\rightarrow\infty$ as $N\rightarrow\infty$, with $t$ growing sufficiently slowly. We will see that such a growth rate is more than enough.
\end{remark}

Write
\[
F=F_t+\Delta_t,
\qquad
F_t=S_{t,t}(F,\boldsymbol{\cdot},\boldsymbol{\cdot}).
\]
Since $F$ is HBV, $F$ is bounded and, by Theorem~\ref{thm:HBV-bounded}, $F_t$ is bounded. Hence, $\Delta_t$ is bounded, with a bound independent of $t$.
\begin{remark}
For convenience we write $F(n,m)$ to mean $F(f(n),f(m))$, where $f$ is
a Steinhaus random multiplicative function. We use the same convention for $F_t$ and
$\Delta_t$.
\end{remark}

We now prove the following estimate.
\begin{lemma}\label{lem:truncation-delta-second-moment}
$\EE\abs{\sum_{n<N}\Delta_t(n,N-n)}^2=o(N)$ as $N\rightarrow\infty$, and so $\abs{\frac{\sum_{n<N}\Delta_t(n,N-n)}{\sqrt{N}}}\xrightarrow{p}0$ as $N\rightarrow\infty$.
\end{lemma}
\begin{proof}
    
We expand the square and separate the diagonal and anti-diagonal contributions \blue{(with an $O_F(1)$ term for the central element, if it exists)}:
\begin{align*}
\EE\abs{\sum_{n<N}\Delta_t(n,N-n)}^2
&=\sum_{n<N}\sum_{m<N}\EE\sqbra{\Delta_t(n,N-n)\conj{\Delta_t(m,N-m)}}\\
&=\sum_{n<N}\EE\abs{\Delta_t(n,N-n)}^2
+\sum_{n<N}\EE\sqbra{\Delta_t(n,N-n)\conj{\Delta_t(N-n,n)}}\\
&+O_F(1)+
\sum_{\substack{n,m<N\\n\neq m, n\neq N-m}}
\EE\sqbra{\Delta_t(n,N-n)\conj{\Delta_t(m,N-m)}}.
\end{align*}

We bound these three terms in turn.
\begin{enumerate}
    \item Let $C_F$ be the uniform bound for $\Delta_t$. Observe that if $n$ and $N-n$ are
independent, then
\[
\EE\abs{\Delta_t(n,N-n)}^2=\EE\abs{\Delta_t(X,Y)}^2=o(1)
\]
as $\Delta_t\ra 0$ in $L^2$ since the Fourier series of a bounded function converges in $L^2$.
We will also use the following counting lemma.
\begin{lem}\label{lem:count-dependent-n}
There are $O_\ep\bra{N^{\frac{1}{3}+\ep}}$ numbers $n<N$ such that\ $n$ and $N-n$ are
dependent.
\end{lem}
\begin{proof}
Suppose $n$ and $N-n$ are dependent. Then, for some integers $\eta_1,\eta_2$ (not
both of them are zero) we have that $n^{\eta_1}(N-n)^{\eta_2}=1$. Moreover, we can
assume $\gcd(\eta_1,\eta_2)=1$.

If one of them is zero, then the other is not, and so either $n=1$ or $N-n=1$, i.e. $O(1)$ solutions here.

If both of them are $\pm1$ then the only possible solution is $n=\frac{N}{2}$, and
so $O(1)$ solutions here.

If one of them is $\pm 2$ and the other is $\pm 1$ then we also have $\ll1$ solutions. 

Suppose now at least one of them is not $\pm1,\pm 2$, and both are non-zero, then as they are
coprime, either $n$ or $N-n$ should be a perfect $k$-th power with $k\geqslant 3$ of a natural number. Then, we
have up to
\[
2\sum_{k=3}^{\log_2N}N^{\frac{1}{k}}=O(\sqrt[3]{N}\log N)=O_\ep\bra{N^{\frac{1}{3}+\ep}}
\]
possible values of $n$ and the lemma follows.
\end{proof}

Therefore,
\[
\sum_{n<N}\EE\abs{\Delta_t(n,N-n)}^2
\leqslant N\EE\abs{\Delta_t(X,Y)}^2+O_\ep\bra{N^{\frac{1}{3}+\ep}}\cdot C_F^2
=o(N).
\]
    
    \item
    \begin{align*}
\abs{\sum_{n<N}\EE\sqbra{\Delta_t(n,N-n)\conj{\Delta_t(N-n,n)}}}
\leqslant&
\frac{1}{2}\sum_{n<N}\EE\sqbra{\abs{\Delta_t(n,N-n)}^2+\abs{\Delta_t(N-n,n)}^2}\\
=&\sum_{n<N}\EE\sqbra{\abs{\Delta_t(n,N-n)}^2}
=o(N),
\end{align*}
    by the previous case.

    \item 
    For the last term
\[
\sum_{\substack{n,m<N\\n\neq m, n\neq N-m}}
\EE\sqbra{\Delta_t(n,N-n)\conj{\Delta_t(m,N-m)}}
\]
note that if $n,m,N-n,N-m$ are not strongly dependent, then by Lemma~\ref{lem:independence-lemma-strong}
we will have $\EE\sqbra{\Delta_t(n,N-n)\conj{\Delta_t(m,N-m)}}=0$.

Indeed, suppose, for example, that $n$ is not in the span of all the others (other
cases are considered similarly). As $F$ has zero marginal expectations, so does
$\Delta_t$, and then
\begin{align*}
&\EE\sqbra{
\Delta_t(n,N-n)\conj{\Delta_t(m,N-m)}
\,\Big|\,
f(N-n),f(m),f(N-m)
} \\
&\qquad=
\EE\sqbra{
\Delta_t(X,f(N-n))\conj{\Delta_t(f(m),f(N-m))}
\,\Big|\,
f(N-n),f(m),f(N-m)
}
= 0.
\end{align*}
Hence, $\EE\sqbra{\Delta_t(n,N-n)\conj{\Delta_t(m,N-m)}}=0$.

    So, this term is bounded by
\[
C_F^2\cdot\#\{(n,m)\text{ s.t. $n\neq m$, $n+m\neq N$ and $n,m,N-n,N-m$ are strongly dependent}\}.
\]
But by Theorem~\ref{thm:second-moment-one}, there are only up to $o(N)$ such pairs,
and so this term is also $o(N)$.
\end{enumerate}
This concludes the proof of the lemma.
\end{proof}

Hence, it is enough to establish
the Central Limit Theorem result for $\sum_{n<N}F_t(n,N-n)$.
\subsubsection{The result for truncated Fourier sum}
Expanding $F_t$ into its Fourier series gives
\[
F_t(f(n),f(N-n))
=\sum_{\abs{k},\abs{m}\leqslant t}c_{k,m}\,f(n)^k\,f(N-n)^m .
\]
Note that the Fourier coefficients $c_{k,m}$ are bounded by the maximum value of $F$.

For a rational non-zero number $r$, let $\Pp(r)$ denote the largest prime number which has
non-zero exponent in $r$; we also set $\Pp(\pm1)=0$.

We adopt the strategy based on McLeish’s theorem that has been used in earlier work on random multiplicative functions, particularly in \cite{Klurman2023} and \cite{chinis2025random}, to establish the Random Chowla conjecture for Steinhaus and Rademacher random multiplicative functions.

For $p$ prime (or $p=0$), define
\[
M_p
=\frac{1}{\sqrt{NC}}\cdot
\sum_{\substack{n<N\\ -t\leqslant k,m\leqslant t\\ \Pp\bra{n^k(N-n)^m}=p}}
c_{k,m}\,f\bra{n^k(N-n)^m}\footnote{We can extend $f$ to rational numbers by multiplicativity},
\]
where $C$ is exactly $\var(F(X,Y))+\cov\bra{F(X,Y), F(Y,X)}$ from
Theorem~\ref{thm:main-steinhaus}. Note that, since $F$ is real-valued and
$\Pp(r)=\Pp\bra{\frac{1}{r}}$, the random variable $M_p$ is real. Moreover,
$(M_p)_{p>0}$ is a martingale difference sequence and
\[
\frac{\sum_{n<N}F_t(n,N-n)}{\sqrt{CN}}
=\sum_{p\geqslant0}M_p(N).
\]

\medskip

\noindent\textbf{The term $M_0$.}
Note that $M_0$ is not stochastic. We claim that $M_0\rightarrow0$ as
$N\rightarrow\infty$, and hence it will be enough to prove the result for the sum
without $M_0$. Indeed,
\begin{align*}
M_0
&=\frac{1}{\sqrt{NC}}\cdot
\sum_{\substack{n<N\\ -t\leqslant k,m\leqslant t\\ n^k(N-n)^m=1}} c_{k,m}\\
&=O_{F}\bra{\frac{1}{\sqrt{N}}\cdot
\sum_{\substack{-t\leqslant k,m\leqslant t\\ k,m\neq 0}}
\#\{n<N\text{ such that }n^k(N-n)^m=1\}},
\end{align*}
where we used the fact that the Fourier coefficients are bounded. Note that, since $k,m\neq0$,
the rational function $n^k(N-n)^m$ cannot be constantly one, and hence we have up to
$\abs{k}+\abs{m}=O(t)$ solutions for each pair $(k,m)$. Using $t=O(\log\log\log N)$, we obtain
\[
M_0
=O_{F}\bra{\frac{1}{\sqrt{N}}\cdot
\sum_{\substack{-t\leqslant k,m\leqslant t\\ k,m\neq 0}}
\#\{n<N\text{ such that }n^k(N-n)^m=1\}}
=O_{F}\bra{\frac{t^3}{\sqrt{N}}}
=o_{F}(1).
\]
From now on, we therefore consider only $\sum_{0<p<N}M_p$.

\medskip

\noindent\textbf{McLeish's theorem.}
As mentioned above, we will use McLeish’s central limit theorem for multi-indexed martingales.
\begin{thm}[McLeish's CLT - {\cite[{Corollary 2.13}]{McLeish}}]
\label{McLeish_CLT}
Let $(k_N)_N\subset \N$ and suppose $(X_{i,N})_i$ is a \textit{martingale difference sequence} on $(\Omega,\mathcal{F},(\mathcal{F}_{i,N})_i,\mathbb{P})$ for $1\leqslant i\leqslant k_N$. Set $S_N:=\sum_{i\leqslant k_N}X_{i,N}$ and suppose the following conditions hold:
\begin{enumerate}
    \item (Normalised Variances) $\sum_{i\leqslant k_N}\EE\left[X_{i,N}^2\right]\rightarrow 1,$ as $N\rightarrow \infty$.
    \item (Lindeberg Condition) For all $\varepsilon>0$, $\sum_{i\leqslant k_N}\EE\left[X_{i,N}^2\mathbbm{1}_{|X_{i,N}|>\varepsilon}\right]\rightarrow 0,$ as $N\rightarrow \infty$.
    \item (Cross-terms Condition) $\limsup_{N\rightarrow \infty} \left(\sum_{i,j\leqslant k_N, i\neq j}\EE\left[X_{i,N}^2X_{j,N}^2\right]\right)\leqslant 1$.
\end{enumerate}
Then: 
$$S_N\xrightarrow{d}\mathcal{N}(0,1),$$
as $N\rightarrow\infty;$ that is, $S_N$ converges in distribution to a normal random variable with mean $0$ and variance $1$.

\end{thm}
As noted in Remark 2.1 in \cite{chinis2025random}, the Lindeberg condition can be replaced by the stronger
\begin{enumerate}
      \setcounter{enumi}{1}
        \item (Fourth-Moment condition) $\sum_{i\leqslant k_N}\EE\left[X_{i,N}^4\right]\rightarrow 0,$ as $N\rightarrow \infty$.
\end{enumerate}

We apply this theorem to the martingale difference sequence $(M_p)_{0<p<N}$.

\medskip
\noindent\textbf{Normalised variances.}
We start from
\begin{align*}
\sum_{0<p<N}\EE\sqbra{\abs{M_p}^2}
=\frac{1}{CN}\cdot
\sum_{\substack{n_1,n_2<N\\ -t\leqslant k_1,k_2,m_1,m_2\leqslant t\\
\Pp(n_1^{k_1}(N-n_1)^{m_1})=\\=\Pp(n_2^{k_2}(N-n_2)^{m_2})\neq0}}
c_{k_1,m_1}\conj{c_{k_2,m_2}}\,
\EE\sqbra{f\bra{\frac{n_1^{k_1}(N-n_1)^{m_1}}{n_2^{k_2}(N-n_2)^{m_2}}}}.
\end{align*}
This implies that
\begin{align*}
\sum_{0<p<N}\EE\sqbra{\abs{M_p}^2}
&=\frac{1}{CN}\cdot
\sum_{\substack{n_1,n_2<N\\ -t\leqslant k_1,k_2,m_1,m_2\leqslant t\\
n_1^{k_1}(N-n_1)^{m_1}=n_2^{k_2}(N-n_2)^{m_2}\neq1}}
c_{k_1,m_1}\conj{c_{k_2,m_2}}\\
&=-M_0^2+\frac{1}{CN}\cdot
\sum_{\substack{n_1,n_2<N\\ -t\leqslant k_1,k_2,m_1,m_2\leqslant t\\
n_1^{k_1}(N-n_1)^{m_1}=n_2^{k_2}(N-n_2)^{m_2}}}
c_{k_1,m_1}\conj{c_{k_2,m_2}}.
\end{align*}

Note that we have a non-zero contribution only in the case when $k_i,m_i\neq 0$, and
in this case, from the proof of Theorem~\ref{thm:second-moment-one} it follows that for each $k_i$ and $m_i$ we can have up to $O_\ep(N^{\frac{1}{2}+\ep})$ non-trivial solutions; since $t=O_\ep(N^\ep)$, the total contribution of non-trivial solutions is $O_{F,\ep}(\frac{1}{N}\cdot t^4\cdot N^{\frac{1}{2}+\ep})=o(1)$. Also, the contribution of $M_0$ is $o(1)$, as shown above.

We obtain trivial solutions either when $n_1=n_2$, $k_1=k_2$ and $m_1=m_2$, or when
$n_1=N-n_2$, $k_1=m_2$ and $k_2=m_1$. Therefore,
\begin{align*}
\sum_{0<p<N}\EE\sqbra{\abs{M_p}^2}
&=o(1)+\frac{1}{CN}\cdot N\cdot
\bra{\sum_{-t\leqslant k_1,m_1\leqslant t}c_{k_1,m_1}\conj{c_{k_1,m_1}}
+\sum_{-t\leqslant k_1,m_1\leqslant t}c_{k_1,m_1}\conj{c_{m_1,k_1}}}\\
&=o(1)+\frac{1}{C}\cdot
\bra{\sum_{\abs{k},\abs{m}\leqslant t}\abs{c_{k,m}}^2
+\sum_{\abs{k},\abs{m}\leqslant t}c_{k,m}\conj{c_{m,k}}}.
\end{align*}

To deduce the normalised variance condition, it is enough to verify that
\[
C=\sum_{k,m}\abs{c_{k,m}}^2+\sum_{k,m}c_{k,m}\conj{c_{m,k}}.
\]
Since $\EE[F(X,y)]=\EE[F(x,Y)]=0$, this follows by Parseval's theorem:
\begin{align*}
\sum_{k,m}\abs{c_{k,m}}^2+\sum_{k,m}c_{k,m}\conj{c_{m,k}}
&=\EE[F(X,Y)^2]+\EE\sqbra{F(X,Y)F(Y,X)}\\
&=\var(F(X,Y))+\cov\bra{F(X,Y), F(Y,X)}\\
&=C.
\end{align*}
Hence, the normalised variance condition holds.

\medskip

\noindent\textbf{Fourth-Moment condition.} 

\begin{align*}\sum_{0<p<N}\EE\bra{\abs{M_p}^4}&=\frac{1}{C^2N^2}\cdot\sum_{p<N}\sum_{\substack{0<n_1,n_2,n_3,n_4<N\\-t\leqslant m_1,\dots,m_4,k_1,\dots,k_4\leqslant t\\\Pp(n_i^{k_i}(N-n_i)^{m_i})=p}}c_{k_1,m_1}c_{k_2,m_2}\conj{c_{k_3,m_3}c_{k_4,m_4}}\cdot\\&\cdot\EE\sqbra{f\bra{n_1^{k_1}(N-n_1)^{m_1}n_2^{k_2}(N-n_2)^{m_2}}\conj{f\bra{n_3^{k_3}(N-n_3)^{m_3}n_4^{k_4}(N-n_4)^{m_4}}}}\end{align*}
The expectation equals $1$ if and only if
\[
n_1^{k_1}(N-n_1)^{m_1}n_2^{k_2}(N-n_2)^{m_2}
=n_3^{k_3}(N-n_3)^{m_3}n_4^{k_4}(N-n_4)^{m_4}
\qquad\tag{$\dagger$}
\]
and equals $0$ otherwise. Hence, we may rewrite
\begin{align*}
\sum_{0<p<N}\EE\bra{\abs{M_p}^4}
=\frac{1}{C^2N^2}\cdot
\sum_{\substack{0<n_1,n_2,n_3,n_4<N\\
-t\leqslant m_1,\dots,m_4,k_1,\dots,k_4\leqslant t\\
\Pp(n_i^{k_i}(N-n_i)^{m_i})=p\text{ for some }p\\
(\dagger)\text{ holds}}}
c_{k_1,m_1}c_{k_2,m_2}\conj{c_{k_3,m_3}c_{k_4,m_4}}.
\end{align*}
Since the Fourier coefficients are bounded, we obtain
\[
\sum_{0<p<N}\EE\bra{\abs{M_p}^4}
=O_{F}\bra{\frac{1}{N^2}\cdot
\sum_{\substack{0<n_1,n_2,n_3,n_4<N\\
-t\leqslant m_1,\dots,m_4,k_1,\dots,k_4\leqslant t\\
m_i,k_i\neq0\\
\Pp(n_i^{k_i}(N-n_i)^{m_i})=p\text{ for some }p\\
(\dagger)\text{ holds}}}1}.
\]

Now apply Theorem~\ref{thm:fourth-moment}. For given $m_i$ and $k_i$, the
contribution from non-trivial solutions is  $\ll_\ep N^{2-\mu+\ep}$ where $\mu>0$ is a constant. Since we
have $O(t^8)=O_\ep(N^\ep)$ different combinations of $m_i$ and $k_i$, the total
contribution of non-trivial solutions is $o(N^2)$. Thus, it is enough to consider
only the solutions where $(\dagger)$ holds trivially.

There are $12$ different combinations where $(\dagger)$ holds trivially (three
different splittings of $n_1,n_2,n_3,n_4$ into pairs, and in each pair two different
equalities). These can be obtained from one another by permuting indices and
replacing $n_i$ by $N-n_i$ for some indices if necessary. Hence, it is enough to
treat one representative case. For ease of notation, consider the case
$n_1=n_3$, $n_2=n_4$ and, since we should have equality as polynomials,
$k_1=k_3$, $m_1=m_3$ and $k_2=k_4$, $m_2=m_4$. The total contribution of this term is
\[
O_F\Biggl(
\frac{1}{N^2}
\sum_{\substack{-t\leqslant k_1,k_2\leqslant t\\-t\leqslant m_1,m_2\leqslant t\\ k_i,m_i\neq 0}}
\#\Bigl\{(n_1,n_2)\in(0,N)^2\cap\mathbb{Z}^2 \Bigm|
\Pp\!\bigl(n_1^{k_1}(N-n_1)^{m_1}\bigr)
=
\Pp\!\bigl(n_2^{k_2}(N-n_2)^{m_2}\bigr)
\neq 0
\Bigr\}
\Biggr).
\]

Fix $k_1,m_1,k_2,m_2$, and bound the number of pairs $(n_1,n_2)$ such that
$\Pp(n_1^{k_1}(N-n_1)^{m_1})=\Pp(n_2^{k_2}(N-n_2)^{m_2})$. We split into two cases.
\begin{itemize}
    \item $\Pp(n_1^{k_1}(N-n_1)^{m_1})=\Pp(n_2^{k_2}(N-n_2)^{m_2})<\sqrt{\log N}$.

    We will denote the contribution of all primes dividing $N$ in a natural number $k$
by $\gcd(k,N^\infty)$ and call it the \emph{$N$-part} of $k$.
    Then $\frac{n_1}{\gcd(n_1,N^\infty)}$ should be $\sqrt{\log N}$-smooth, and by formula $(1.18)$ from \cite{granville08}, $\frac{n_1}{\gcd(n_1,N^\infty)}$ has up to
$\bra{\frac{\log N}{\sqrt{\log N}}}^{(1+o(1))\pi\bra{\sqrt{\log N}}}$ different possibilities, and we will show in Claim~\ref{claim:bound-radicals} below, that $\gcd(n_1,N^\infty)$ has $\ll_\ep N^\ep$ possibilities. Observe that
\begin{align*}
\bra{\frac{\log N}{\sqrt{\log N}}}^{(1+o(1))\pi\bra{\sqrt{\log N}}}\cdot O_\ep(N^\ep)
&=\bra{\log N}^{\frac{1+o(1)}{2}\cdot \frac{\sqrt{\log N}}{\frac{1}{2}\log\log N}}\cdot O_\ep(N^\ep)=\\
&=\exp\bra{(1+o(1))\sqrt{\log N}}\cdot O_\ep(N^\ep)
=O_\ep(N^\ep).
\end{align*}
Similarly for $n_2$, so we have up to $O_\ep(N^\ep)$ pairs here.

    \item $\Pp(n_1^{k_1}(N-n_1)^{m_1})=\Pp(n_2^{k_2}(N-n_2)^{m_2})=p\geqslant\sqrt{\log N}$.

    Then either $n_1$ or $N-n_1$ is divisible by $p$, so for a given $p$ there are up to $\frac{2N}{p}$ possibilities for $n_1$. Similarly for $n_2$, so the total number of pairs is bounded by \[
\sum_{p\geqslant\sqrt{\log N}}\frac{4N^2}{p^2}\leqslant \sum_{n\geqslant\sqrt{\log N}}\frac{4N^2}{n^2}
=O\bra{\frac{N^2}{\sqrt{\log N}}}.
\]
\end{itemize}

Note that we have $O(t^4)=o(\sqrt{\log N})$ different combinations of
$k_1,m_1,k_2,m_2$. Therefore, in total, we have $o(N^2)$ trivial terms in the sum for one combination,
and summing over all $12$ combinations yields
\[
\sum_{0<p<N}\EE\bra{\abs{M_p}^4}
=\frac{1}{C^2N^2}\cdot o(N^2)
=o(1),
\]
which proves the Fourth-Moment condition.

\medskip

\noindent\textbf{Cross-term condition.} Since we already know that $\sum_{0<p<N}\EE\bra{\abs{M_p}^4}=o(1)$, the cross-term
condition is equivalent to
\[
\limsup_{N\rightarrow\infty}\EE\sqbra{\bra{\sum_{0<p<N}M_p^2}^2}\leqslant 1.
\]

We can rewrite
\begin{align*}
\EE\sqbra{\bra{\sum_{0<p<N}M_p^2}^2}=&\frac{1}{C^2N^2}\cdot\sum_{\substack{0<n_1,n_2,n_3,n_4<N\\-t\leqslant m_1,\dots,m_4,k_1,\dots,k_4\leqslant t\\\Pp(n_1^{k_1}(N-n_1)^{m_1})=\Pp(n_2^{k_2}(N-n_2)^{m_2})>0\\\Pp(n_3^{k_3}(N-n_3)^{m_3})=\Pp(n_4^{k_4}(N-n_4)^{m_4})>0}}c_{k_1,m_1}c_{k_2,m_2}c_{k_3,m_3}c_{k_4,m_4}\cdot\\&\cdot\EE\sqbra{f\bra{n_1^{k_1}(N-n_1)^{m_1}n_2^{k_2}(N-n_2)^{m_2}n_3^{k_3}(N-n_3)^{m_3}n_4^{k_4}(N-n_4)^{m_4}}}\end{align*}
The expectation equals $1$ if
\[
n_1^{k_1}(N-n_1)^{m_1}n_2^{k_2}(N-n_2)^{m_2}
n_3^{k_3}(N-n_3)^{m_3}n_4^{k_4}(N-n_4)^{m_4}=1
\qquad\tag{$\ddagger$}
\]
and equals $0$ otherwise. Hence, the sum becomes
\[
\frac{1}{C^2N^2}\cdot
\sum_{\substack{0<n_1,n_2,n_3,n_4<N\\
-t\leqslant m_1,\dots,m_4,k_1,\dots,k_4\leqslant t\\
\Pp(n_1^{k_1}(N-n_1)^{m_1})=\Pp(n_2^{k_2}(N-n_2)^{m_2})\\
\Pp(n_3^{k_3}(N-n_3)^{m_3})=\Pp(n_4^{k_4}(N-n_4)^{m_4})\\
(\ddagger)\text{ holds}}}
c_{k_1,m_1}c_{k_2,m_2}c_{k_3,m_3}c_{k_4,m_4}.
\]

As in the fourth-moment argument, it is enough to consider only the trivial
solutions of $(\ddagger)$. Moreover, in our splitting of $n_1,n_2,n_3,n_4$ into pairs we need
$n_1$ and $n_2$ to lie in the same pair; otherwise we obtain
\[
\Pp(n_1^{k_1}(N-n_1)^{m_1})
=\Pp(n_2^{k_2}(N-n_2)^{m_2})
=\Pp(n_3^{k_3}(N-n_3)^{m_3})
=\Pp(n_4^{k_4}(N-n_4)^{m_4}),
\]
and from the proof of the fourth-moment condition, we already know that the
contribution of such solutions is \textit{small}. Thus, we have only four contributing types of
(trivial) solutions:
\begin{itemize}
        \item
When $k_1=-k_2$, $m_1=-m_2$ and $k_3=-k_4$, $m_3=-m_4$. The contribution is
\[
\frac{1}{C^2}\cdot\bra{\sum_{\abs{k},\abs{m}\leqslant t}\abs{c_{k,m}}^2}^2.
\]
\item
When $k_1=-k_2$, $m_1=-m_2$ and $k_3=-m_4$, $m_3=-k_4$. The contribution is
\[
\frac{1}{C^2}\cdot
\bra{\sum_{\abs{k},\abs{m}\leqslant t}\abs{c_{k,m}}^2}\cdot\,
\bra{\sum_{\abs{k},\abs{m}\leqslant t}c_{k,m}\conj{c_{m,k}}}.
\]
    \item
When $k_1=-m_2$, $m_1=-k_2$ and $k_3=-k_4$, $m_3=-m_4$. The contribution is
\[
\frac{1}{C^2}\cdot
\bra{\sum_{\abs{k},\abs{m}\leqslant t}c_{k,m}\conj{c_{m,k}}}\cdot\,
\bra{\sum_{\abs{k},\abs{m}\leqslant t}\abs{c_{k,m}}^2}.
\]
    \item When $k_1=-m_2$, $m_1=-k_2$ and $k_3=-m_4$, $m_3=-k_4$. The contribution is
\[
\frac{1}{C^2}\cdot
\bra{\sum_{\abs{k},\abs{m}\leqslant t}c_{k,m}\conj{c_{m,k}}}^2.
\]
\end{itemize}
Summing these contributions gives
\begin{align*}
\EE\sqbra{\bra{\sum_{0<p<N}M_p^2}^2}
&=o(1)+
\bra{\frac{
\sum_{\abs{k},\abs{m}\leqslant t}\abs{c_{k,m}}^2
+\sum_{\abs{k},\abs{m}\leqslant t}c_{k,m}\conj{c_{m,k}}
}{C}}^2=1+o(1),
\end{align*}
since we showed above that
$C=\sum_{k,m}\abs{c_{k,m}}^2+\sum_{k,m}c_{k,m}\conj{c_{m,k}}$. Therefore,
\[
\limsup_{N\ra\infty}\EE\sqbra{\bra{\sum_{0<p<N}M_p^2}^2}=1,
\]
and the cross-term condition holds.

\medskip

\noindent\textbf{Conclusion.} Having verified all three conditions of Theorem~\ref{McLeish_CLT}, we can conclude that
$\sum_{p>0} M_p(N)$ tends to $\mathcal{N}(0,1)$ in distribution. Since $M_0=o(1)$, it
follows that
\[
\frac{\sum_{n<N}F_t(n,N-n)}{\sqrt{CN}}
\xrightarrow{d}\mathcal{N}(0,1),
\]
which establishes the result for the truncated Fourier series. Combining with the results from Section~\ref{subsection:truncation}, we conclude that Theorem~\ref{thm:main-steinhaus} holds in the case $G=F$, that is when $c_{n,0}=c_{0,m}=0$ for all $n,m\in\Z$.
\subsubsection{Deducing  Theorem~\ref{thm:main-steinhaus} in the most general case}\label{subsection:deducing-steinhaus-main}
Define
\[
u(x):=\EE F(x,Y)-\EE F(X,Y),
\qquad
v(y):=\EE F(X,y)-\EE F(X,Y).
\]
Then
\begin{equation}
F(f(n),f(N-n))
=
G(f(n),f(N-n))
+u(f(n))
+v(f(N-n))
+\EE F(X,Y).
\end{equation}
In Theorem~\ref{thm:main-steinhaus} we explicitly remove the contribution of the
constant term. Note that if $F$ has HBV, then $G$ has HBV as well, so the
convergence result for the $G$-term has already been established above. Therefore,
to conclude the proof of Theorem~\ref{thm:main-steinhaus}, it remains to show that
\[
\sum_{n<N} u(f(n))
\quad\text{and}\quad
\sum_{n<N} v(f(N-n))=\sum_{n<N} v(f(n))
\]
exhibit better than square-root cancellation.
\begin{remark}
If $C$ from Theorem~\ref{thm:main-steinhaus} is zero, then
$G(X,Y)=-G(Y,X)$ almost everywhere. Thus, after pairing the terms indexed by
$n$ and $N-n$, all terms in
\[
\sum_{n<N} G(f(n), f(N-n))
\]
cancel except possibly those for which $n$ and $N-n$ are dependent. By
Lemma~\ref{lem:count-dependent-n}, there are
$O_\ep(N^{\frac{1}{3}+\ep})$ such indices. Since $G$ is bounded, we obtain
\[
\sum_{n<N} G(f(n), f(N-n))
=O_{G,\ep}(N^{\frac{1}{3}+\ep}).
\]
Taking $\ep<\frac{1}{6}$, we indeed obtain better than square-root cancellation.
\end{remark}
Thus, it is enough to prove the following theorem, which will finish the proof of Theorem~\ref{thm:main-steinhaus} (subject to the second and fourth moment theorems established in the subsequent sections).
\begin{thm}\label{thm:one-variable-cancellation}
Let $u$ be a measurable real-valued function on $\T$ with bounded harmonic
variation and mean zero. Then, if $f$ is a Steinhaus random multiplicative
function, we have better than square-root cancellation for sums $\sum_{n<N}u(f(n))$, i.e.
\[
\frac{\sum_{n<N}u(f(n))}{\sqrt{N}}\xrightarrow[]{p}0.
\]
\end{thm}
\begin{proof}
    We use Theorem~1 from \cite{harper2020moments}.
    \begin{thm}[Harper]\label{thm:harper}
Let $f$ be a Steinhaus random multiplicative function. Then,
\[
\EE\abs{\sum_{n\leqslant x} f(n)}
\asymp
\frac{\sqrt{x}}{\sqrt[4]{\log\log x}}.
\]
\end{thm}
\begin{remark}
    \blue{In fact, we only need that the expectation grows more slowly than $\sqrt{x}$.}
\end{remark}

Let $(d_k)_{k\in\Z}$ be the Fourier coefficients of $u$. Since $u$ has mean zero, we
have $d_0=0$. Let $t=\lfloor\log\log\log N\rfloor$ and let $u_t$ be the truncation of the Fourier
series keeping the indices between $-t$ and $t$. Write
\[
u=u_t+\delta_t.
\]
Then $\delta_t\to 0$ in $L^2$ as $u$ is bounded. Moreover, both $u_t$ and $\delta_t$ are real-valued.

\medskip
\noindent\textbf{Step 1: the truncated part $u_t$.}
Since $f^k$ (for $k\neq 0$) is also a Steinhaus random multiplicative function, Theorem~\ref{thm:harper} implies
\[
\EE\abs{\sum_{n<N}f(n)^k}
=O\!\bra{\frac{\sqrt{N}}{\sqrt[4]{\log\log N}}}
\qquad (k\neq 0).
\]
Using the Fourier expansion of $u_t$ and the triangle inequality,
\begin{align*}
\EE\abs{\sum_{n<N}u_t(f(n))}
&\leqslant
\sum_{\abs{k}\leqslant t}\abs{d_k}\cdot
\EE\abs{\sum_{n<N}f(n)^k}.
\end{align*}
Since the coefficients $d_k$ are bounded (as $u$ is bounded) and $d_0=0$, we obtain
\begin{align*}
\EE\abs{\sum_{n<N}u_t(f(n))}
&=O\!\bra{\frac{t\sqrt{N}}{\sqrt[4]{\log\log N}}}
=O\!\bra{\frac{\log\log\log N\cdot\sqrt{N}}{\sqrt[4]{\log\log N}}}
=o(\sqrt{N}).
\end{align*}
In particular, by Markov's inequality,
\[
\frac{\sum_{n<N}u_t(f(n))}{\sqrt{N}}\xrightarrow[]{p}0.
\]

\medskip
\noindent\textbf{Step 2: the error term $\delta_t$.}
Consider
\begin{align*}
\EE\bra{\sum_{n<N}\delta_t(f(n))}^2
&=\EE\bra{\sum_{n,m<N}\delta_t(f(n))\delta_t(f(m))}\\
&=\delta_t(1)^2
+(N-2)\EE\delta_t(X)^2
+\EE\bra{\sum_{\substack{1<n,m<N\\ n\neq m}}\delta_t(f(n))\delta_t(f(m))},
\end{align*}
where $X$ is a uniform random point on $\T$. Since $u$ has bounded harmonic variation, the functions $\delta_t$ are uniformly bounded; let $C_u$ be a uniform bound for $\delta_t$.

As before, if $n$ and $m$ are not dependent, then the corresponding expectation is
zero. We claim that there are only $o(N)$ dependent pairs with $1<n,m<N$ and
$n\neq m$.

Indeed, if $n$ and $m$ are dependent, then there exist coprime $a,b>0$ such that
$n^a=m^b$. Hence, $n=s^b$ and $m=s^a$ for some $1<s<N$. Note that $a,b\leqslant\log_2N$.
Moreover, $s\leqslant N^{\min\{\frac{1}{a},\frac{1}{b}\}}$. Since $n\neq m$, we cannot
have $a=b=1$, and therefore $s\leqslant \sqrt{N}$. It follows that there are
$\ll_\ep N^{\frac{1}{2}+\ep}$ dependent non-equal pairs with $1<n,m<N$.

Since $\delta_t\xrightarrow{L^2}0$, we have $\EE\delta_t(X)^2=o(1)$, and therefore
\begin{align*}
\EE\bra{\sum_{n<N}\delta_t(f(n))}^2
&=O(1)+N\cdot o(1)+C_u^2\cdot O_\ep(N^{\frac{1}{2}+\ep})
=o(N).
\end{align*}
By Markov's inequality,
\[
\frac{\sum_{n<N}\delta_t(f(n))}{\sqrt{N}}\xrightarrow{p}0
\qquad\text{as }N\rightarrow\infty.
\]

\medskip
\noindent\textbf{Step 3: conclusion.}
Combining the bounds for $u_t$ and $\delta_t$ yields
\[
\frac{\sum_{n<N}u(f(n))}{\sqrt{N}}\xrightarrow{p}0
\qquad\text{as }N\rightarrow\infty,
\]
as claimed.
\end{proof}

Note that 
\begin{align*}HV(u)&=\sup_{E\in \Omega}\sum_{[a_i,b_i]\in E} \frac{\abs{\EE F(b_i,Y)-\EE F(a_i,Y)}}{i}\leqslant\\&\leqslant \sup_{E\in \Omega}\sum_{[a_i,b_i]\in E} \frac{\EE \abs{F(b_i,Y)-F(a_i,Y)}}{i}\leqslant HV_1(F)<\infty,\end{align*} so $u$, and similarly, $v$ have bounded harmonic variation, and mean zero by definition. Hence, we may apply
Theorem~\ref{thm:one-variable-cancellation} to both $u$ and $v$, which completes the
argument.
\subsection{Moment estimates}
In this section, we prove several results concerning the computation of the second and fourth moments. We caution that the proofs involve extensive calculations. More precisely, we prove the following three theorems:
\begin{thm}\label{thm:second-moment-one}
There are up to $O_\ep(N^{\frac{1}{2}+\ep})$ pairs $(n,m)$ up to $N$ such that
$n\neq m$, $n+m\neq N$, and $n,m,N-n,N-m$ are strongly dependent.
\end{thm}
\begin{remark}
    In fact, we will also show that for every relation
$n^\alpha(N-n)^\beta m^\gamma(N-m)^\delta=1$ where
$\alpha,\beta,\gamma,\delta\neq 0$, we have
$O_\ep(N^{\frac{1}{2}+\ep})$ solutions, apart from those on
a diagonal $n=m$ or $n+m=N$ along which the relation holds
identically as an equality of rational functions.
\end{remark}
\begin{thm}\label{thm:fourth-moment}
There exists a constant $\mu>0$ such that the following holds. Let
$\alpha_1$, $\alpha_2$, $\beta_1$, $\beta_2$, $\gamma_1$, $\gamma_2$, $\delta_1$, $\delta_2$ be non-zero
integers whose absolute values are less than $\log\log\log N$, and $a,b,c,d\in\{1,2,...,N-1\}$. Then the equation
\begin{equation}\tag{$\ast$}\label{eq:fourth-moment-star}
\begin{aligned}
&a^{\alpha_1}(N-a)^{\alpha_2}b^{\beta_1}(N-b)^{\beta_2}\cdot c^{\gamma_1}(N-c)^{\gamma_2}d^{\delta_1}(N-d)^{\delta_2}=1
\end{aligned}
\end{equation}
has $O_\ep\bra{N^{2-\mu+\ep}}$ non-trivial solutions. We call a solution
\emph{trivial} if $a,b,c,d$ split into pairs such that, in each pair, the numbers are equal to or add up to $N$, and moreover, assuming these conditions are true,
\eqref{eq:fourth-moment-star} holds as an equality of rational functions.
\end{thm}
\begin{remark}
In fact, the $\mu$ obtained here is $0.001666\ldots$, but the analysis presented below can be optimised to improve this constant.
\end{remark}
\begin{thm}\label{thm:second-moment-2}
Suppose $A$ and $B$ are coprime natural numbers such that $AB>1$, $A$ and $B$ are up
to $N^{4\log\log\log N}$ and $N$-smooth. Then, for given non-zero integers
$\alpha,\beta,\gamma,\delta$ up to $\log\log\log N$ in absolute value there are up to
$O_\ep(N^{\frac{1}{2}+\ep})$ solutions to
\[
An^{\alpha}(N-n)^\beta m^\gamma(N-m)^\delta=B,
\qquad n,m<N
\]
uniformly in $A$ and $B$.
\end{thm}
\begin{remark}
The difference between this theorem and Theorem~\ref{thm:second-moment-one} is that
here we assume that we have a non-cancelled contribution from $a$ and $d$.
\end{remark}
\begin{remark}
These conditions on $A$ and $B$ come from the fact that we apply Theorem~\ref{thm:second-moment-2} in the case when
\[
\frac{A}{B}=a^{\alpha_1}(N-a)^{\alpha_2}d^{\delta_1}(N-d)^{\delta_2}.
\]
If we take $A$ and $B$ to be the numerator and denominator of this expression, then
they clearly satisfy all the conditions.
\end{remark}
\subsubsection{Second moment calculation}
To prove Theorem~\ref{thm:second-moment-one}, we prove the following three claims; the theorem will then follow immediately.
\begin{claim}\label{claim:exist-small-relation}
Suppose $n,N-n,m,N-m$ are strongly dependent. Then there exist non-zero integers
$\alpha,\beta,\gamma,\delta$, each of size at most
$O((\log N)^3)=O_\ep(N^\ep)$, such that
\[
n^\alpha (N-n)^\beta m^\gamma (N-m)^\delta = 1.
\]
\end{claim}
\begin{proof}
    We first show that there exist non-zero integers $\alpha,\beta,\gamma,\delta$ with
\[
n^\alpha (N-n)^\beta m^\gamma (N-m)^\delta = 1.
\]
In fact, it suffices to find such $\alpha,\beta,\gamma,\delta\in\Q$.

    Let $V\subset\Q^4$ be the set of all $(\alpha,\beta,\gamma,\delta)$ satisfying the
above relation. Then $V$ is a $\Q$-linear subspace of $\Q^4$. We claim that $V$ is
not contained in the union of the four coordinate hyperplanes
$\{\alpha=0\}\cup\{\beta=0\}\cup\{\gamma=0\}\cup\{\delta=0\}$.
Since $\Q$ is infinite, this containment would be possible only if $V$ were
contained in one of the hyperplanes, say $\{\alpha=0\}$ without loss of generality.
But then $n$ would not lie in the span of $N-n,m,N-m$, contradicting the assumption
that $n,N-n,m,N-m$ are strongly dependent. Thus, there exists a vector in $V$ with
all four coordinates non-zero.

Now let $d$ be the dimension of the span of $n,N-n,m,N-m$. Then $1\leqslant d\leqslant 3$. The corresponding vectors of prime valuations have coefficients bounded by
$O(\log N)$. Hence, by the Bombieri--Vaaler theorem
(Theorem 2 in \cite{bombieri1983}),
the space $V$ admits an integral basis whose coordinates are bounded by $O((\log N)^3)$.

Let $\mathbf{u}_1,\dots,\mathbf{u}_{4-d}$ be such a basis of $V$, and consider
\[
\mathbf{v}=\lambda_1\mathbf{u}_1+\cdots+\lambda_{4-d}\mathbf{u}_{4-d},
\qquad \lambda_i\in\Z.
\]
Each coordinate of $\mathbf{v}$ is a linear polynomial in $(\lambda_i)_i$, so the
product of the coordinates is a polynomial
\[
P(\lambda_1,\dots,\lambda_{4-d})
\]
of degree at most $4$. Since $V$ contains a vector with all coordinates non-zero,
this polynomial is not identically zero. Therefore, by the Combinatorial
Nullstellensatz, there is a choice of $\lambda_i\in\{0,1,2,3,4\}$ such that
$P(\lambda_1,\dots,\lambda_{4-d})\neq 0$. For this choice, $\mathbf{v}$ has all
coordinates non-zero and bounded by $O((\log N)^3)$, as required.
\end{proof}
\begin{claim}\label{claim:bound-radicals}
\[
\#\{k<N:\rad k\mid \rad m\}=O_\ep(N^\ep)
\qquad\text{for } m=N^{O(1)}.
\]
\end{claim}
\begin{proof}
    Suppose $m<N^A$ for some fixed $A$. Write
\[
m=q_1^{\alpha_1}\cdots q_r^{\alpha_r},
\qquad
q_1<q_2<\cdots<q_r,
\qquad
\alpha_i>0.
\]

    It suffices to prove the claim in the case where the prime factors are the first
$r$ primes, i.e.\ $q_i=p_i$ for $i=1,\dots,r$. Indeed, consider any integer $t$
such that $\rad t\mid \rad m$. Then we may write
\[
t=q_1^{\beta_1}\cdots q_r^{\beta_r}.
\]
Define $t':=p_1^{\beta_1}\cdots p_r^{\beta_r}$. Since $p_i\leqslant q_i$ for all $i$,
we have $t'\leqslant t$. Consequently, if $t<N$ then $t'<N$. Moreover,
$\rad t'\mid \rad m'$, where
\[
m':=p_1^{\alpha_1}\cdots p_r^{\alpha_r}\leqslant m.
\]
Since the correspondence $t\mapsto t'$ is one-to-one and $m'\leqslant m<N^A$, it is
enough to establish the desired bound for $m'$, i.e.\ we may assume $q_i=p_i$.
    
    Thus, suppose $m=2^{\alpha_1}\cdots p_r^{\alpha_r}$. Then $\rad t\mid \rad m$ if and only if $t$ is $p_r$-smooth. Let $\Psi(x,y)$ denote the number of $y$-smooth integers up to $x$.

    Since $m\leqslant N^A$, we have $2^r\leqslant N^A$, hence
$r\leqslant A\log_2N<2A\log N$, and therefore $p_r=O_A(\log N\log\log N)$. We now
apply formula (1.19) from \cite{granville08}:
\begin{align*}
\log\Psi(N,p_r)
&\leqslant \log\Psi(N,C_A\log N\log\log N)\\
&=\frac{\log N}{\log \bra{C_A\log N\log\log N}}
  \, g\!\bra{\frac{C_A\log N\log\log N}{\log N}}\\
&\qquad \times
  \bra{1+O\!\bra{
    \frac{1}{\log\bra{C_A\log N\log\log N}}
    +\frac{1}{\log\log N}}},\\
&\qquad \text{where } g(\kappa)=\log(1+\kappa)+\kappa\log\!\bra{1+\frac{1}{\kappa}}.
\end{align*}
    In particular,
\[
\log\Psi(N,p_r)
=O\bra{\frac{\log N}{\log\log N}\cdot\log\log\log N\cdot 1}
=o(\log N),
\]
and hence $\Psi(N,p_r)=O_\ep(N^\ep)$, which proves the claim.
\end{proof}
\begin{remark}
The claim above, in fact, shows that there are $O_\ep(N^\ep)$ different values of
$\gcd(k,N^\infty)$ for $k<N$.
\end{remark}
\begin{claim}\label{claim:fixed-parts-fixed-exponents}
For given non-zero $\alpha,\beta,\gamma,\delta$ up to $O\bra{(\log N)^3}$, and
$a,b,c,d<N$, there are only $O_\ep\bra{N^{\frac{1}{2}+\ep}}$ pairs $(n,m)$ satisfying
the conditions of Theorem~\ref{thm:second-moment-one}, and moreover
\[
\gcd(n,N^\infty)=a,\quad
\gcd(N-n,N^\infty)=b,\quad
\gcd(m,N^\infty)=c,\quad
\gcd(N-m,N^\infty)=d,
\]
and
\[
n^\alpha (N-n)^\beta m^\gamma (N-m)^\delta=1
\qquad (\dagger).
\]
\end{claim}
\begin{remark}
Assuming the three claims, we obtain at most
\[
O_\ep(N^{\ep})\cdot O_\ep(N^{\ep})\cdot O_\ep\bra{N^{\frac{1}{2}+\ep}}
=O_\ep\bra{N^{\frac{1}{2}+\ep}}
\]
pairs $(n,m)$ satisfying the conditions of Theorem~\ref{thm:second-moment-one}.
\end{remark}
\begin{proof}
The proof proceeds by a case analysis. First, we may assume that
$\gcd(\alpha,\beta,\gamma,\delta)=1$. Note that $\alpha,\beta,\gamma,\delta$ cannot
all have the same sign, as this would force $(\dagger)$ to have no solutions.
Consequently, the signs are distributed either as a $3{:}1$ split or as a $2{:}2$
split. In the latter case, we must further distinguish whether the coefficients of
the same sign correspond to the same variable or not.
    \begin{enumerate}
        \item\textbf{A $3{:}1$ split.}
Without loss of generality, assume $\alpha>0$ and $\beta,\gamma,\delta<0$. Moving
terms across in $(\dagger)$, we may equivalently solve, with all
$\alpha,\beta,\gamma,\delta>0$,
\[
n^\alpha=(N-n)^\beta m^\gamma (N-m)^\delta.
\]
Then
\[
\gcd(n^\alpha,N^\infty)=\gcd((N-n)^\beta m^\gamma (N-m)^\delta,N^\infty),
\]
which implies $a^\alpha=b^\beta c^\gamma d^\delta$, and hence
\[
\bra{\frac{n}{a}}^\alpha=
\bra{\frac{N-n}{b}}^\beta
\bra{\frac{m}{c}}^\gamma
\bra{\frac{N-m}{d}}^\delta.
\] 
        
Note that if a prime $p$ divides both $n$ and $N-n$, then $p\mid N$. Hence
$\frac{n}{a}$ and $\frac{N-n}{b}$ are coprime. As $\beta>0$ this implies
$\frac{N-n}{b}=1$, so there is only one possible value of $n$. Moreover,
$m^\gamma\mid n^\alpha$, so $\rad m\mid \rad n$. By Claim~\ref{claim:bound-radicals},
for a given $n$ there are $O_\ep(N^\ep)$ possible values of $m$. Thus, there are
$O_\ep(N^\ep)$ solutions in this case.
        \item \textbf{A $2{:}2$ split, different variables on each side.}
Without loss of generality, suppose (after rearranging) that all
$\alpha,\beta,\gamma,\delta>0$ and
\[
n^\alpha m^\gamma=(N-n)^\beta(N-m)^\delta.
\]
        As before, $a^\alpha c^\gamma=b^\beta d^\delta$, and we also have both
$\bra{\frac{n}{a}}^\alpha\mid\bra{\frac{N-m}{d}}^\delta$ and
$\bra{\frac{N-m}{d}}^\delta\mid\bra{\frac{n}{a}}^\alpha$, which is only possible if
\[
\bra{\frac{n}{a}}^\alpha=\bra{\frac{N-m}{d}}^\delta.
\]
Similarly,
\[
\bra{\frac{N-n}{b}}^\beta=\bra{\frac{m}{c}}^\gamma.
\]
In particular, if one of $n$ and $m$ is known, then the other has at most one
possible value.

We consider three subcases.
\begin{itemize}
\item If $\alpha\neq\delta$.

Assume $\alpha<\delta$. Then $\frac{n}{a}$ is a perfect
$\frac{\delta}{\gcd(\alpha,\delta)}\geqslant 2$ power, and therefore $\frac{n}{a}$
admits at most
\[
O\bra{N^{\frac{\gcd(\alpha,\delta)}{\delta}}}
=O\bra{N^{\frac{1}{2}}}
\]
possible values. The case $\alpha>\delta$ is analogous.

\item If $\beta\neq\gamma$, we argue as in the previous case.

\item If $\alpha=\delta$ and $\beta=\gamma$.

Then
\[
\begin{cases}
\frac{N-n}{b}=\frac{m}{c}\\[2pt]
\frac{N-m}d=\frac{n}{a}
\end{cases}
\Longleftrightarrow
\begin{cases}
\frac{n}{b}+\frac{m}{c}=\frac{N}{b}\\[2pt]
\frac{n}{a}+\frac{m}{d}=\frac{N}{d}
\end{cases}
\]
\begin{itemize}
\item If $bd\neq ac$, the system has at most one solution, so the desired bound
holds.

\item If $bd=ac$ and the system has at least one solution, then
$\frac{N}{b}\cdot\frac{b}{a}=\frac{N}{d}$, so $a=d$, and similarly $b=c$. This
implies $n+m=N$, but we count only non-trivial solutions.
\end{itemize}
\end{itemize}
    \item \textbf{A $2{:}2$ split, same variable on each side.}
Without loss of generality (after rearranging), assume all
$\alpha,\beta,\gamma,\delta>0$ and
\[
n^\alpha (N-n)^\beta=m^\gamma(N-m)^\delta.
\]

As before, we deduce
\[
\bra{\frac{n}{a}}^\alpha\bra{\frac{N-n}{b}}^\beta
=
\bra{\frac{m}{c}}^\gamma\bra{\frac{N-m}{d}}^\delta,
\quad
\gcd\bra{\frac{n}{a},\frac{N-n}{b}}
=
\gcd\bra{\frac{m}{c},\frac{N-m}{d}}=1.
\]
Also note that $\rad m\mid \rad\bigl(n(N-n)\bigr)$, so by
Claim~\ref{claim:bound-radicals}, knowing $n$ bounds sufficiently well the number of
possible values of $m$, and vice-versa. 

Define
\[
\begin{aligned}
x &= \gcd\bra{\bra{\frac{n}{a}}^\alpha,\bra{\frac{m}{c}}^\gamma},
\qquad &
z &= \gcd\bra{\bra{\frac{n}{a}}^\alpha,\bra{\frac{N-m}{d}}^\delta},\\
y &= \gcd\bra{\bra{\frac{N-n}{b}}^\beta,\bra{\frac{m}{c}}^\gamma},
\qquad &
w &= \gcd\bra{\bra{\frac{N-n}{b}}^\beta,\bra{\frac{N-m}{d}}^\delta}.
\end{aligned}
\]
Then $x,y,z,w$ are pairwise coprime and
\[
\begin{aligned}
\bra{\frac{n}{a}}^\alpha&=xz,
\qquad &
\bra{\frac{m}{c}}^\gamma&=xy,\\
\bra{\frac{N-n}{b}}^\beta&=yw,
\qquad &
\bra{\frac{N-m}{d}}^\delta&=zw.
\end{aligned}
\]
Moreover,
\[
\begin{aligned}
x\text{ and }z\text{ are natural }\alpha&\text{-th powers},
\qquad &
x\text{ and }y\text{ are natural }\gamma&\text{-th powers},\\
y\text{ and }w\text{ are natural }\beta&\text{-th powers},
\qquad &
z\text{ and }w\text{ are natural }\delta&\text{-th powers}.
\end{aligned}
\]

Suppose $\gamma\nmid\alpha$ and $\delta\nmid\alpha$. Since $x$ is both an $\alpha$-th
and $\gamma$-th power,
\[
x=\bra{x'}^{\frac{\alpha\gamma}{\gcd(\alpha,\gamma)}},
\]
and similarly,
\[
z=\bra{z'}^{\frac{\alpha\delta}{\gcd(\alpha,\delta)}}.
\]
Therefore,
\[
\frac{n}{a}
=\bra{x'}^{\frac{\gamma}{\gcd(\alpha,\gamma)}}
 \bra{z'}^{\frac{\delta}{\gcd(\alpha,\delta)}}
\geqslant (x'z')^2,
\]
since $\gamma\nmid\alpha$ and $\delta\nmid\alpha$. Hence, $x'z'\leqslant\sqrt{N}$, and
the number of possible pairs $(x',z')$ is
\[
O\bra{N^{\frac{1}{2}}\log N}
=O_\ep\bra{N^{\frac{1}{2}+\ep}}.
\]
Each such pair yields a unique value of $n$, so we obtain
$O_\ep\bra{N^{\frac{1}{2}+\ep}}$ solutions overall.

Thus, we may assume $\gamma\mid\alpha$ or $\delta\mid\alpha$. Similarly, we obtain
the following relations (the first one is included for reference):
\[
\left[
\begin{array}{ll}
\gamma\mid\alpha\\
\delta\mid\alpha
\end{array}\right.
\qquad
\left[
\begin{array}{ll}
\gamma\mid\beta\\
\delta\mid\beta
\end{array}\right.
\qquad
\left[
\begin{array}{ll}
\alpha\mid\gamma\\
\beta\mid\gamma
\end{array}\right.
\qquad
\left[
\begin{array}{ll}
\alpha\mid\delta\\
\beta\mid\delta
\end{array}\right.
\]
We now use the following lemma.

\begin{lemma}\label{lem:divisibility-structure}
These relations imply that one of the following is true:
$\alpha=\gamma$, $\alpha=\delta$, $\beta=\gamma$, or $\beta=\delta$.
\end{lemma}
\begin{proof}[Proof of the lemma]
Without loss of generality, assume $\gamma\mid \alpha$, and none of the equalities
above holds. Then $\alpha\nmid\gamma$, hence $\beta\mid\gamma$. Since $\beta\neq\gamma$,
we have $\gamma\nmid\beta$, so $\delta\mid\beta$. Similarly, $\beta\nmid\delta$, hence
$\alpha\mid\delta$. Overall we obtain
\[
\gamma\mid\alpha\mid\delta\mid\beta\mid\gamma,
\]
so $\alpha=\beta=\gamma=\delta$, which is impossible.
\end{proof}
Assume without loss of generality that $\alpha=\gamma$. We first treat the case of small exponents.
\begin{lemma}\label{lem:small-exponents}
Suppose $\alpha,\beta,\gamma,\delta\leqslant \log\log N$. Then there are
$O_\ep\bra{N^{\frac{1}{2}+\ep}}$ solutions.
\end{lemma}
\begin{proof}
In this case, we can use the Bombieri--Pila theorem.

\begin{thm}[Bombieri--Pila, Theorem 5 from \cite{Bombieri}]\label{thm:bombieri-pila}
Let $\mathcal{C}$ be an absolutely irreducible curve of degree
$r\geqslant 2$. If $N\geqslant \exp(r^6)$, then the number of integral points on
$\mathcal{C}$ inside the square $[0,N]\times[0,N]$ does not exceed
\[
N^{\frac{1}{r}}\exp\bigl(12\sqrt{r\log N\log \log N}\bigr).
\]
\end{thm}

    Consider
\[
P(n,m)=n^\alpha(N-n)^\beta-m^\gamma(N-m)^\delta.
\]
It suffices to bound the number of integral points on $P(n,m)=0$ in $[0,N]^2$
excluding, if necessary, the trivial families $n=m$ and $n+m=N$. Let $Q$ be a non-linear
irreducible factor of $P$ over $\C$. Then
\[
r=\deg Q\leqslant \deg P=O(\log\log N),
\]
so in particular $N\geqslant\exp(r^6)$. Hence, we may apply Bombieri--Pila to $Q$,
obtaining
\[
N^{\frac{1}{r}}\exp\bigl(12\sqrt{r\log N\log \log N}\bigr)
=O_\ep\bigl(N^{\frac{1}{2}+\ep}\bigr),
\]
solutions of $Q(n,m)=0$ in $[0,N]^2$ since $2\leqslant r\ll\log\log N$.

Now suppose $Q$ is a linear factor of $P$. We may assume $Q(X,Y)=X-pY-q$ (if the
coefficient of $X$ is zero, we swap $n$ and $m$). Then $P(pY+q,Y)=0$ for all $Y$,
i.e.
\[
(pY+q)^\alpha(N-pY-q)^\beta=Y^\gamma(N-Y)^\delta.
\]
Therefore, either $Y$ is proportional to $pY+q$ and $N-Y$ is proportional to
$N-pY-q$, or $Y$ is proportional to $N-pY-q$ and $N-Y$ is proportional to $pY+q$.
In the first case, we get $q=0$ and $pN=N$, i.e.\ $Q=X-Y$. In the second case, we get
$N-q=0$ and $q=-pN$, i.e.\ $Q=X+Y-N$. In both cases, we obtain only the trivial
solutions.
\begin{remark}
If $X-Y$ is a factor of $P$, then $\alpha=\gamma$ and $\beta=\delta$. If $X+Y-N$ is a factor of $P$, then $\alpha=\delta$ and $\beta=\gamma$. This confirms that we cannot obtain ``\textit{many}'' trivial solutions for a non-obvious choice of $\alpha,\beta,\gamma,\delta$.
\end{remark}

Thus, the number of non-trivial integral points on $P$ is bounded by
\[
\sum_{\substack{Q\mid P\\ \deg Q\geqslant 2}}
O_\ep\bra{N^{\frac{1}{2}+\ep}}
=O_\ep\bra{\deg P\cdot N^{\frac{1}{2}+\ep}}
=O_\ep\bra{\log\log N\cdot N^{\frac{1}{2}+\ep}}
=O_\ep\bra{N^{\frac{1}{2}+\ep}}.
\]
\end{proof}

Assume now that the hypothesis of Lemma~\ref{lem:small-exponents} is false, so that $\max\{\alpha,\beta,\gamma,\delta\}>\log\log N$. We first dispose of the case $\beta\neq\delta$. If $\alpha\nmid\delta$ then $\beta\mid\delta$. If $\delta\mid\beta$ then
$\alpha=\gamma$ and $\beta=\delta$ and we are done. Otherwise, $\gamma\mid\beta\mid\delta$
and so $\alpha\mid\delta$ since $\alpha=\gamma$, a contradiction. Now assume $\alpha\mid\delta$ and $\beta\neq\delta$. Either $\gamma\mid\beta$ or
$\delta\mid\beta$. But $\alpha$ divides both $\gamma$ and $\delta$, so
$\alpha\mid\gcd(\alpha,\beta,\gamma,\delta)$. Since we assumed this gcd is $1$, we
must have $\alpha=\gamma=1$.

We may assume without loss of generality that $\beta<\delta$ (the roles of $\beta$
and $\delta$ are symmetric here), so $\delta\geqslant\log\log N$.

We know that $z$ and $w$ are perfect $\delta$-th powers, and $y$ and $w$ are perfect $\beta$-th powers. Therefore, $w$ is a perfect $\frac{\beta\delta}{\gcd(\beta,\delta)}$-th power. Write
\[
z=\bra{z'}^\delta,
\qquad
w=\bra{w'}^{\frac{\beta\delta}{\gcd(\beta,\delta)}}.
\]
Since $\bra{z'}^\delta\leqslant \frac{n}{a}<N$, the number of possibilities for $z$
is $O\bra{N^{\frac{1}{\delta}}}=O_\ep(N^\ep)$ as $\delta\geqslant \log\log N$.

Also,
\[
w^{\frac{1}{\beta}}\leqslant\frac{N-n}{b}<N,
\qquad\text{i.e.}\qquad
\bra{w'}^{\frac{\delta}{\gcd(\beta,\delta)}}<N.
\]
Since $\beta<\delta$, this implies $\frac{\delta}{\gcd(\beta,\delta)}\geqslant 2$ and $\bra{w'}^2<N$, so $w$ has 
$\ll N^{\frac{1}{2}}$ possibilities. Hence, there are
$O_\ep\bra{N^{\frac{1}{2}+\ep}}$ possible pairs $(z,w)$, and each such pair
determines $m$ uniquely. This yields the required bound in this case.

Therefore, the only remaining case is
\[
\alpha=\gamma,\qquad \beta=\delta,\qquad
\max\{\alpha,\beta,\gamma,\delta\}>\log\log N.
\]
Moreover, $\gcd(\alpha,\beta)=1$, and we may assume without loss of generality $\alpha\leqslant\beta$, so
$\beta\geqslant\log\log N$.

As before, $y$ and $z$ are perfect $\alpha\beta$-th powers. Since
$y\leqslant N^\gamma=N^\alpha$, the number of possibilities for $y$ is
$O\bra{N^{\frac{1}{\beta}}}=O_\ep(N^\ep)$, and similarly for $z$. Thus, there are
$O_\ep(N^\ep)$ possible pairs $(y,z)$. For a given pair $(y,z)$,
\begin{align*}
\bra{\frac{n}{m}\cdot\frac{c}{a}}^\alpha&=\frac{z}{y},\\
\bra{\frac{N-n}{N-m}\cdot\frac{d}{b}}^\beta&=\frac{y}{z},
\end{align*} so the pair $(y,z)$ specifies $\frac{n}{m}$ and $\frac{N-n}{N-m}$ uniquely.
Unless $\frac{n}{m}=\frac{N-n}{N-m}=1$ (which yields only trivial solutions), the
system has at most one solution. Indeed,
\[
\begin{cases}
\frac{n}{m}=A,\\[2pt]
\frac{N-n}{N-m}=B
\end{cases}
\Longleftrightarrow
\begin{cases}
n-Am=0,\\
n-Bm=N-NB.
\end{cases}
\]
This system can have more than one solution only if $A=B$ and $N-NB=0$, i.e.\ if $A=B=1$. Therefore, in this remaining case, we also have 
$\ll_\ep N^{\frac{1}{2}+\ep}$ non-trivial solutions, completing the proof.
\end{enumerate}
\end{proof}
\subsubsection{Fourth moment calculation}

\begin{proof}[Proof of Theorem~\ref{thm:fourth-moment}]
    Call $a^{\alpha_1}(N-a)^{\alpha_2}$ the \emph{contribution of $a$}. Define similarly the contributions of $b,c$, and $d$.

    \medskip
    \noindent\textbf{Step 1: the cancelling-two-contributions case.}
    
    First, suppose that the product of the contributions for two variables is equal to $1$. This case covers all trivial solutions. We show that there are $\ll_\ep N^{\frac{3}{2}+\ep}$ \emph{non-trivial} solutions in this scenario. Without loss of generality, assume these two variables are $a$ and $d$. Then we have
the system
    \begin{align*}
a^{\alpha_1}(N-a)^{\alpha_2}d^{\delta_1}(N-d)^{\delta_2}&=1,\\
b^{\beta_1}(N-b)^{\beta_2}c^{\gamma_1}(N-c)^{\gamma_2}&=1.
\end{align*}

    By Theorem~\ref{thm:second-moment-one}, each of these equations has
$O_\ep\bra{N^{\frac{1}{2}+\ep}}$ non-trivial solutions (and possibly $O(N)$ trivial
ones). Since we do not count solutions in which \emph{both} pairs are trivial, we
obtain  $\ll_\ep N^{\frac{3}{2}+\ep}$ solutions here.

    Thus, from now on, we assume that the contribution from no pair of variables cancels
out; in particular, we no longer need to account for trivial solutions.

\medskip
\noindent\textbf{Step 2: a lemma on small radicals.}
We will use the following bound.
    \begin{lem}\label{lem:bound-radicals}
Suppose $\zeta$ and $C>0$ are fixed. Then,
\[
\#\{x\leqslant M:\text{\emph{rad}}(x)\leqslant CM^\zeta\}=O_{\zeta,C,\ep}(M^{\zeta+\ep}).
\]
\end{lem}
    \begin{proof}
        The statement is trivial for $\zeta=0$ and for $\zeta\geqslant 1$. So assume $0<\zeta<1$. 
        
We use the formulation given in \cite{Robert}. Introduce the notation
\[
N(x,y):=\#\{n\leqslant x\text{ such that }\rad n\leqslant y\}.
\]
Combining formulas \emph{(1.12)} and \emph{(2.12)} from \cite{Robert}, we obtain the following.
\begin{lem}
Let \(x\geqslant 16\) and assume that
\[
e^{(\log x)^{3/4+\varepsilon}} \leqslant y \leqslant x .
\]
Set \(v=\log(x/y)\). Then
\[
N(x,y)=yF(v)\bra{1+O\bra{\sqrt{\frac{\log(v+2)}{\log x}}}},
\]
where \(F\) is as in \emph{(1.9)} of \cite{Robert}. Moreover, for \(v\geqslant 3\),
\[
F(v)=\exp\!\left(\sqrt{\frac{8v}{\log v}}\,(1+o(1))\right)
\qquad (v\to\infty).
\]
\end{lem}
Substitute $x=M$ and $y=CM^\zeta$. Since $0<\zeta<1$, we have $e^{(\log x)^{3/4+\varepsilon}} \leqslant y \leqslant x $, $v\asymp \log M$ and
$\frac{\log(v+2)}{\log x}=O(1)$. Therefore,
\begin{align*}
N(M,CM^\zeta)
&=CM^\zeta\exp\bra{O\bra{\sqrt{\frac{\log M}{\log \log M}}}}\cdot O(1)\\
&=O(M^\zeta)\cdot O_\ep(M^\ep)
=O_{\zeta,C,\ep}\bra{M^{\zeta+\ep}},
\end{align*}
as required.
\end{proof}

\medskip
\noindent\textbf{Step 3: ordering by the largest radical and a case split.}
Assume without loss of generality that $\rad(d)\rad(N-d)$ is maximal among the four quantities $\rad(x)\rad(N-x)$, $x\in\{a,b,c,d\}$.

Fix parameters $0\leqslant\eta<\xi\leqslant 2$ (to be chosen later) and consider
solutions with
\[
N^\eta\leqslant \rad(d)\rad(N-d)\leqslant N^\xi.
\]
Also fix $0<\lambda<\nu<1$. Along with the proof, we record constraints on
$\lambda,\nu,\eta,\xi$; we mark each such constraint by $(\star)$.

We consider three cases.
    \begin{itemize}
    \item \textbf{Case 1:} $\rad\gcd(d,N)=\rad\gcd(N-d,N)<N^\lambda$.

Note that either $\rad(d)\leqslant N^{\frac{\xi}2}$ or
$\rad(N-d)\leqslant N^{\frac{\xi}2}$. Hence, by Lemma~\ref{lem:bound-radicals}, $d$
has $\ll_\ep N^{\frac{\xi}2+\ep}$ possible values.

Let $p$ be a prime such that $p\mid d(N-d)$ but $p\nmid N$. Then $p$ cannot divide
both $d$ and $N-d$, and since $\delta_1,\delta_2\neq 0$, the contribution of $p$ in
the contribution of $d$ is not $1$. Hence, $p$ must divide one of the remaining
contributions. Considering all such primes, we obtain
\[
X=\frac{\rad(d)\rad(N-d)}{\bra{\rad(\gcd(d,N))}^2}\mid\rad\bra{a(N-a)b(N-b)c(N-c)}.
\]
Note that $X$ is coprime to $N$ and $X\geqslant N^{\eta-2\lambda}$.
There exists a divisor $Y$ of $X$ such that $Y\mid \rad(a(N-a))$ or $\rad(b(N-b))$ or $\rad(c(N-c))$, and moreover, $Y\geqslant \sqrt[3]{X}$, so in particular 
$Y\geqslant N^{\frac{\eta-2\lambda}{3}}$. Without loss of generality,
assume $Y\mid \rad(a(N-a))$.

We need the following bound.

\begin{lem}\label{lem:Y-divides-aNa}
Let $Y$ be coprime to $N$ and $Y\leqslant N^{C}$ for some fixed $C>0$. Then the
number of $a$ such that $Y\mid a(N-a)$ is
$
\ll_\ep\bra{N^{\ep}\cdot\max\{1;\frac{N}{Y}\}}.
$
\end{lem}

\begin{proof}
Note that
\[
\#\{a\mod Y: a(N-a)\equiv 0\mod Y\}= 2^{\omega(Y)},
\]
since both sides are multiplicative and, for $Y$ a prime power coprime to $N$, both
are equal to $2$. Also, since $Y\leqslant N^C$, we have $2^{\omega(Y)}=O_\ep(N^\ep)$.
Hence, the number of such $a$ is bounded by
\[
O_\ep(N^\ep)\cdot\bra{1+\frac{N}{Y}},
\]
which implies the claimed bound.
\end{proof}
    Applying Lemma~\ref{lem:Y-divides-aNa} and noting that
$\frac{\eta-2\lambda}{3}\leqslant \frac{2}{3}$, we obtain that the number of possible
values of $a$, given $d$, is
\[
O_\ep\bra{N^{1-\frac{\eta-2\lambda}{3}+\ep}}.
\]
\begin{remark}
    Note that we have $O_\ep(N^\ep)$ possible values of $Y$ which does not affect the bound.
\end{remark}
Finally, by Theorem~\ref{thm:second-moment-2}, which will be proven in the following section, there are
$O_\ep(N^{\frac{1}{2}+\ep})$ possible pairs $(b,c)$ given $a$ and $d$. Altogether
this gives
\[
O_\ep\bra{N^{\frac{3}{2}+\frac{\xi}2-\frac{\eta-2\lambda}{3}+\ep}}
\]
solutions in Case~1. We gain over the exponent $2$ provided that
\[
\frac{3}{2}+\frac{\xi}2-\frac{\eta-2\lambda}{3}<2.\qquad\tag{$\star$}
\]
\item \textbf{Case 2:} $N^\lambda\leqslant \rad\gcd(d,N)\leqslant N^\nu$.

Let $d'=\frac{d}{\gcd(d,N^\infty)}$ and $d''=\frac{N-d}{\gcd(N-d,N^\infty)}$. Then 
\begin{align*}
\rad\bra{d'}\cdot\rad\bra{d''}=\frac{\rad(d)\rad(N-d)}{(\rad\gcd(d,N))^2}\leqslant N^{\xi-2\lambda}
=\bra{N^{1-\lambda}}^{\frac{\xi-2\lambda}{1-\lambda}}.
\end{align*}
Note that $d',d''\leqslant\frac{N}{\gcd(d,N)}\leqslant N^{1-\lambda}$. So, assuming
\[
\xi-2\lambda>0\qquad\tag{$\star$}
\]
and repeating the same argument as in Case~1, we find that either $d'$ or $d''$ has
\[
O_\ep\bra{\bra{N^{1-\lambda}}^{\frac{\xi-2\lambda}{2(1-\lambda)}+\ep}}
=O_\ep\bra{N^{\frac{\xi}2-\lambda+\ep}}
\]
possible values. Then by Claim~\ref{claim:bound-radicals}, both $\gcd(d,N^\infty)$ and $\gcd(N-d,N^\infty)$ have $\ll_\ep N^\ep$ possible values. It follows that $d$ has $O_\ep\bra{N^{\frac{\xi}2-\lambda+\ep}}$ possibilities in this case.

Similarly to Case~1, we obtain that $a$ has
$O_\ep\bra{N^{1-\frac{\eta-2\nu}3+\ep}}$ possibilities, and $(b,c)$ has
$O_\ep\bra{N^{\frac{1}{2}+\ep}}$ possibilities given $(a,d)$. Hence, the total number   
of solutions in Case~2 is
\[
O_\ep\bra{N^{\frac{3}{2}+\frac{\xi}2-\lambda-\frac{\eta-2\nu}{3}+\ep}}.
\]
We gain over $2$ provided that
\[
\frac{3}{2}+\frac{\xi}2-\lambda-\frac{\eta-2\nu}{3}<2.\qquad\tag{$\star$}
\]
\item \textbf{Case 3:} $\rad\gcd(d,N)>N^\nu$.

Arguing as above, and assuming
\[
\xi-2\nu>0\qquad\tag{$\star$},
\]
we obtain $O_\ep\bra{N^{\frac{\xi}2-\nu+\ep}}$ possibilities for $d$ in this case.
Since $\rad(d)\rad(N-d)$ is maximal, we have $\rad(a)\rad(N-a)\leqslant N^\xi$, and
by the same argument as in Case~1, $a$ has $O_\ep\bra{N^{\frac{\xi}{2}+\ep}}$
possibilities. Given $(a,d)$, there are, by Theorem~\ref{thm:second-moment-2},
$O_\ep\bra{N^{\frac{1}{2}+\ep}}$ possibilities for the pair $(b,c)$. Therefore, in
Case~3 we obtain
\[
O_\ep\bra{N^{\xi-\nu+\frac{1}{2}+\ep}}
\]
solutions. We gain over $2$ provided that
\[
\xi-\nu+\frac{1}{2}<2.\qquad\tag{$\star$}
\]
    \end{itemize}

\medskip
\noindent\textbf{Step 4: parameter choices.}
Summarising, if $\lambda,\nu,\eta,\xi$ satisfy
\[
\begin{cases}
0\leqslant\eta<\xi\leqslant2,\\
0<\lambda<\nu<1,\\
\xi>2\lambda,\\
\xi>2\nu,\\
\frac{3}{2}+\frac{\xi}2-\frac{\eta-2\lambda}{3}<2,\\
\frac{3}{2}+\frac{\xi}2-\lambda-\frac{\eta-2\nu}{3}<2,\\
\xi-\nu+\frac{1}{2}<2,
\end{cases}
\]
then there exists $\mu(\eta,\xi,\lambda,\nu)>0$ such that the number of solutions
with
\[
N^\eta\leqslant\rad(d)\rad(N-d)\leqslant N^\xi
\]
is $O_\ep\bra{N^{2-\mu(\eta,\xi,\lambda,\nu)+\ep}}$.

    Consider the following choices (one checks directly that the above constraints hold
in each case):
\begin{itemize}
\item $\xi=2,\eta=1.95$ and $\lambda=0.22,\nu=0.55$,
\item $\xi=1.95,\eta=1.8$ and $\lambda=0.18,\nu=0.455$,
\item $\xi=1.8,\eta=1.45$ and $\lambda=0.12,\nu=0.302$.
\end{itemize}
Hence, there is some $\mu>0$ such that the number of solutions with
$\rad(d)\cdot \rad(N-d)\geqslant N^{1.45}$ is $O_\ep(N^{2-\mu+\ep})$.

Finally, if $\rad(d)\cdot \rad(N-d)\leqslant N^{1.45}$ then, by the same argument as
in Case~1, $d$ has $O_\ep(N^{0.725+\ep})$ possibilities. Also
$\rad(a)\cdot \rad(N-a)\leqslant N^{1.45}$, so $a$ has
$O_\ep(N^{0.725+\ep})$ possibilities. Given $(a,d)$, there are
$O_\ep(N^{0.5+\ep})$ possible pairs $(b,c)$ by Theorem~\ref{thm:second-moment-2}. Thus, in this range we have
$\ll_\ep N^{1.95+\ep}$ solutions. Updating $\mu$ if necessary completes the proof.
\end{proof}
\subsubsection{Second moment calculation II}
\begin{proof}[Proof of Theorem~\ref{thm:second-moment-2}]

    As in Theorem~\ref{thm:second-moment-one}, we may assume that
\[
a=\gcd(n,N^\infty),\qquad b=\gcd(N-n,N^\infty),\qquad
c=\gcd(m,N^\infty),\qquad d=\gcd(N-m,N^\infty)
\]
are fixed. As before, we argue by considering the possible sign patterns of
$\alpha,\beta,\gamma,\delta$. Note that here we also allow the case in which all
four exponents have the same sign.
    \begin{itemize}
        \item \textbf{All of $\alpha,\beta,\gamma,\delta$ have the same sign.}
Without loss of generality, assume they are positive.

Then $n\mid B$. Recall that
\[
\log\tau(B)=O\bra{\frac{\log B}{\log \log B}}.
\]
Hence, as $B\leqslant N^{4\log\log\log N}$,
\[
\log\tau(B)
=O\bra{\frac{\log N\cdot \log\log\log N}{\log \log N}}
=o(\log N).
\]
Therefore, $\tau(B)=O_\ep(N^\ep)$, and hence $n$ has $O_\ep(N^\ep)$ possibilities.
Similarly for $m$, so there are $O_\ep(N^\ep)$ solutions in this case.

\item \textbf{A $3{:}1$ split.}
Without loss of generality, we may assume we are solving
\[
An^\alpha=B(N-n)^\beta m^\gamma(N-m)^\delta,
\qquad \alpha,\beta,\gamma,\delta>0.
\]
Then $\frac{N-n}{b}\mid A$, and by an argument similar to the one in the previous case, $n$ has
$O_\ep(N^\ep)$ possibilities. Having fixed $n$, we can apply the same argument to
$m$, and we obtain $O_\ep(N^\ep)$ solutions in this case.

\item \textbf{A $2{:}2$ split, different variables on each side.}
Without loss of generality, assume we are solving
\[
An^\alpha m^\gamma=B(N-n)^\beta(N-m)^\delta,
\qquad \alpha,\beta,\gamma,\delta>0,
\]
and moreover assume $A>1$.

We have
\begin{align*}
\bra{\frac{n}{a}}^\alpha &\mid B\bra{\frac{N-m}{d}}^\delta,\\
\bra{\frac{N-m}{d}}^\delta &\mid A\bra{\frac{n}{a}}^\alpha.
\end{align*}
Since $A$ and $B$ both have $O_\ep(N^\ep)$ divisors, it follows that there are
$O_\ep(N^\ep)$ possible values of $\frac{n^\alpha}{(N-m)^\delta}$, and similarly
$O_\ep(N^\ep)$ possible values of $\frac{(N-n)^\beta}{m^\gamma}$.

Fix
\[
x=\frac{n^\alpha}{(N-m)^\delta},
\qquad
y=\frac{(N-n)^\beta}{m^\gamma}.
\]
Then
\[
N-m=x^{-\frac{1}{\delta}}n^{\frac{\alpha}{\delta}},
\qquad
m=y^{-\frac{1}{\gamma}}(N-n)^{\frac{\beta}{\gamma}},
\]
and hence
\[
x^{-\frac{1}{\delta}}n^{\frac{\alpha}{\delta}}
+y^{-\frac{1}{\gamma}}(N-n)^{\frac{\beta}{\gamma}}
= N.
\]
Equivalently,
\[
\bra{N-x^{-\frac{1}{\delta}}\bra{n^{\frac{1}{\delta}}}^\alpha}^\gamma
=y^{-1}\bra{N-\bra{n^{\frac{1}{\delta}}}^\delta}^\beta.
\]
This is a polynomial equation in $n^{\frac{1}{\delta}}$. If the corresponding
polynomial is not identically zero, then it has
$O(\alpha\gamma+\beta\delta)=O_\ep(N^\ep)$ roots, so $n$ has
$O_\ep(N^\ep)$ possibilities; moreover, using $y=\frac{(N-n)^\beta}{m^\gamma}$ we
have up to $1$ solution for $m$.

Thus, in the non-degenerate subcase, we obtain $O_\ep(N^\ep)$ solutions. Now suppose
that the polynomial
\[
\bra{N-x^{-\frac{1}{\delta}}z^\alpha}^\gamma=y^{-1}\bra{N-z^\delta}^\beta
\]
is identically zero. Taking $z=N^{\frac{1}{\delta}}$, the right-hand side vanishes,
so $N=x^{-\frac{1}{\delta}}N^{\frac{\alpha}{\delta}}$, i.e.
$x=N^{\alpha-\delta}$. This implies
\[
\bra{\frac{n}{N}}^\alpha=\bra{\frac{N-m}{N}}^\delta.
\]
Similarly, taking $z=0$, we obtain that
\[
\bra{\frac{N-n}{N}}^\beta=\bra{\frac{m}{N}}^\gamma.
\]

If $\alpha=\delta$ and $\beta=\gamma$, this implies that $n+m=N$, and from our
equation we deduce that $A=B$. Since $\gcd(A,B)=1$, this gives $A=B=1$, contradicting
$AB>1$.

Hence, either $\alpha\neq\delta$ or $\beta\neq\gamma$, and then one of
\[
\frac{n}{N},\ \frac{N-n}N,\ \frac{m}{N},\ \frac{N-m}{N}
\]
must be a perfect power of a rational number. Without loss of generality, assume
this is $\frac{n}{N}$. Then, $\frac{n}{\gcd(n, N)}$ is a perfect power of a natural number, and so it has $\ll_\ep N^{\frac{1}{2}+\ep}$ possibilities, and then, by the divisor bound, there are $\ll_\ep N^{\frac{1}{2}+\ep} $ possibilities for
$n$. Given $n$, there are $O_\ep(N^\ep)$ possibilities for $m$. Therefore, the
total number of solutions in this case is
$O_\ep(N^{\frac{1}{2}+\ep})$, as required.
\item \textbf{A $2{:}2$ split, same variable on each side.}
Without loss of generality, assume we are solving
\[
An^\alpha (N-n)^\beta=Bm^\gamma(N-m)^\delta,
\qquad \alpha,\beta,\gamma,\delta>0.
\]

View this as a polynomial equation in $n,m$ of degree $O(\log\log\log N)$; denote
the corresponding polynomial by $P$. Note that the only possible linear factors of
$P$ are $n-m$ and $n+m-N$. We claim that the conditions on $A$ and $B$ forbid these
to be factors of $P$. Indeed, if one of them were a factor, then for all $n$ one of
the following identities would hold:
\begin{align*}
An^\alpha(N-n)^\beta &=Bn^\gamma(N-n)^\delta,\\
An^\alpha(N-n)^\beta &=B(N-n)^\gamma n^\delta.
\end{align*}

In either case we would have $A=B$, which is impossible since $\gcd(A,B)=1$ and
$AB>1$.

Since $\deg P=O(\log\log\log N)$, we may apply the Bombieri--Pila theorem and repeat
the analysis from Lemma~\ref{lem:small-exponents}, obtaining
$O_\ep(N^{\frac{1}{2}+\ep})$ solutions.
    \end{itemize}
    In all cases, we obtain the desired bound, completing the proof.
\end{proof}
\section{A note about Goldbach's conjecture for Rademacher RMF}\label{section:goldbach-rademacher}

Since a Rademacher random multiplicative function takes only finitely many values,
every function of the form $F(f(n),f(N-n))$ can be written as a finite
linear combination of indicator functions
\[
\mathbbm{1}(f(n)=c_1)\,\mathbbm{1}(f(N-n)=c_2),
\qquad c_1,c_2\in\{0,\pm1\}.
\]
Thus, in principle, it would suffice to establish an analogue of
Corollary~\ref{cor:product-indicators}. Unfortunately, the method that worked well
for the Steinhaus random multiplicative function does not seem to apply in full
generality in the Rademacher setting. \blue{Nevertheless, as stated in the introduction, we can still prove an analogue of Corollary~\ref{cor:complex-gaussian}, namely Theorem~\ref{thm:main-rademacher}, from which we deduce a partial result for indicator functions, stated as Corollary~\ref{cor:rademacher-indicators}.}

\begin{proof}[Proof of Corollary~\ref{cor:rademacher-indicators}]
Without loss of generality, assume $c_1=1$ and $c_2=-1$. Then
\[
\mathbbm{1}(f(n)=c_1)=\frac{1}{2}\mu^2(n)(f(n)+1),
\qquad
\mathbbm{1}(f(N-n)=c_2)=\frac{1}{2}\mu^2(N-n)(1-f(N-n)).
\]
Therefore,
\begin{align*}
&\frac{\sum_{n<N}\mathbbm{1}(f(n)=c_1)\mathbbm{1}(f(N-n)=c_2)-\frac{S(N)}{4}}
{\sqrt{\frac{S(N)}8}}=\\
&=
\frac{\sum_{n<N}\mu^2(n)\mu^2(N-n)\bra{f(n)-f(N-n)-f(n)f(N-n)}}
{4\cdot\sqrt{\frac{S(N)}8}}\\
&=
\frac{\sum_{n<N}\mu^2(N-n)f(n)-\sum_{n<N}\mu^2(n)f(N-n)-\sum_{n<N}f(n)f(N-n)}
{\sqrt{2S(N)}}.
\end{align*}
The first two terms cancel (this is where we need $c_1c_2=-1$), and we are left with minus
the expression from Theorem~\ref{thm:main-rademacher}.
\end{proof}
\begin{remark}
Note that
\[
\sum_{n<N} f(n)f(N-n)
=2\sum_{n<\frac{N}{2}}f(n)f(N-n)+O(1),
\]
where the $O(1)$ error comes from the possible $f(\frac{N}{2})^2$ term. Also,
\[
S(N)
=2\#\left\{n<\frac{N}2\ \Bigm|\ n\text{ and }N-n\text{ are square-free}\right\}+O(1).
\]
Therefore, Theorem~\ref{thm:main-rademacher} is equivalent to
\[
\frac{\sum_{n<\frac{N}{2}}f(n)f(N-n)}
{\sqrt{\#\{n<\frac{N}2\mid n\text{ and }N-n\text{ are square-free}\}}}
\ \xrightarrow{\,d\,}\ \mathcal{N}(0,1).
\]
We denote $T(N):=\#\Bigl\{n<\frac{N}2\ \Bigm|\ n\text{ and }N-n\text{ are square-free}\Bigr\}$
and observe that $T(N)\asymp N$.
\end{remark}
\subsection{Proof of Theorem~\ref{thm:main-rademacher}}

We follow the same approach as in the proof of Theorem~\ref{thm:main-steinhaus},
using McLeish's theorem (Theorem~\ref{McLeish_CLT}). To implement this method in the
Rademacher setting, we need analogues of
Theorem~\ref{thm:fourth-moment} and Theorem~\ref{thm:second-moment-one}.
These are provided by the following two counting lemmata.
\begin{lem}\label{lem:rademacher-second-moment}
There are $O_\ep(N^\ep)$ pairs $(a,b)$ of distinct positive integers less than $\frac{N}{2}$
such that $a$, $b$, $N-a$, $N-b$ are all square-free, $a(N-a)b(N-b)$ is a perfect
square.
\end{lem}

\begin{lem}\label{lem:rademacher-fourth-moment}
There are $3T(N)^2+O_\ep(N^{2-\frac{4}{11}+\ep})$ quadruples $(a,b,c,d)$ of
positive integers less than $\frac{N}{2}$ such that $a$, $b$, $c$, $d$, $N-a$,
$N-b$, $N-c$, $N-d$ are all square-free and
$a(N-a)b(N-b)c(N-c)d(N-d)$ is a perfect square.
\end{lem}
We will also use the following auxiliary lemma.

\begin{lem}\label{lem:pell-solutions}
Suppose $r$ is a positive rational number, not equal to $1$, and whose
numerator and denominator are $N^{O(1)}$. Then there are $O_\ep(N^\ep)$ pairs of
positive integers $n,m$ less than $N$ such that
\[
\frac{n(N-n)}{m(N-m)}=r.
\]
\end{lem}
\begin{proof}
We consider two cases.

\begin{itemize}
\item \textbf{Case 1: $r$ is not a perfect square.}

We use the following proposition from \cite{Cilleruelo2011}.

\begin{prop}[Integral points on quadratic curves, Proposition $1$ from \cite{Cilleruelo2011}]
\label{prop:Pell}
Let $\abs{A},\abs{B},\abs{C},\abs{D},\abs{E},\abs{F} \leqslant N^{O(1)}$ and assume that $B^2 - 4AC \neq \square$. Then, the Diophantine equation
\[
Ax^2 + Bxy + Cy^2 + Dx + Ey + F = 0
\]
has at most $N^{o(1)}$ integral solutions with $1 \leqslant \abs{x}, \abs{y} \leqslant N^{O(1)}$.
\end{prop}
Write $r=\frac{p}{q}$ where $p$ and $q$ are positive coprime integers. Then
\[
\frac{n(N-n)}{m(N-m)}=r
\Longleftrightarrow
pm^2-qn^2+qN\cdot n-pN\cdot m=0
\tag{$\dagger$}.
\]
The coefficients in $(\dagger)$ are $N^{O(1)}$, and
\[
-4\cdot p\cdot (-q)=4pq=(2q)^2r
\]
is not a square. Therefore, by Proposition~\ref{prop:Pell}, there are
$O_\ep(N^\ep)$ solutions.

\item \textbf{Case 2: $r$ is a perfect square.}

In this case we can write $r=\frac{p^2}{q^2}$ where $p$ and $q$ are positive
coprime integers. Therefore,
\[
\frac{n(N-n)}{m(N-m)}=r
\Longleftrightarrow
(pm-qn)(pm+qn)=N(p^2m-q^2n)
\tag{$\ddagger$}.
\]

Let $x=pm-qn$ and $y=pm+qn$, so $pm=\frac{x+y}{2}$ and $qn=\frac{y-x}{2}$. Note that
$y>0$. Also, $x\neq 0$: otherwise the LHS of $(\ddagger)$ vanishes, and so $pm-qn=p^2m-q^2n=0$, which implies $r=1$, a
contradiction. Rewriting $(\ddagger)$ in terms of $x$ and $y$, we obtain
\[
2xy=N(p(x+y)-q(y-x))
\Longleftrightarrow
2xy=N((p+q)x+(p-q)y).
\]
Hence,
\[
\begin{cases}
x\mid N(p-q)y,\\
y\mid N(p+q)x.
\end{cases}
\]
Let $d=\gcd(x,y)$ and write $x=dx'$, $y=dy'$. Then $x'$ and $y'$ are coprime, and
\[
\begin{cases}
x'\mid N(p-q),\\
y'\mid N(p+q).
\end{cases}
\]
Since $p+q\neq 0$ and $p-q\neq 0$ and they are $N^{O(1)}$, the divisor bound implies
that there are $O_\ep(N^\ep)$ possibilities for the pair $(x',y')$.

Finally, note that
\[
2xy=N(p(x+y)-q(y-x))
\Rightarrow
2dx'y'=N(p(x'+y')-q(y'-x')).
\]
Since $xy\neq 0$, for each admissible choice of $(x',y')$ there is at most one
possible value of $d$. Therefore, there are $O_\ep(N^\ep)$ possibilities for
$(x,y)$, and hence $O_\ep(N^\ep)$ solutions $(n,m)$ to the original equation.
\end{itemize}
\end{proof}
\begin{proof}[Proof of Lemma~\ref{lem:rademacher-second-moment}]
Let
\[
A=\frac{a(N-a)}{\gcd(a,N)^2},
\qquad
B=\frac{b(N-b)}{\gcd(b,N)^2}.
\]
Observe that $A$ and $B$ are both square-free and $AB$ is a perfect square.
Therefore, $A=B$, and hence
\[
\frac{a(N-a)}{b(N-b)}=\frac{\gcd(a,N)^2}{\gcd(b,N)^2}.
\]
Since $N$ has $O_\ep(N^\ep)$ divisors, the ratio $\frac{a(N-a)}{b(N-b)}$ has
$O_\ep(N^\ep)$ possible values, and each value satisfies the bounds of
Lemma~\ref{lem:pell-solutions}.

If $\frac{a(N-a)}{b(N-b)}=1$, then $a=b$. The case $a+b=N$ is impossible since
$a,b<\frac{N}{2}$. In any case, we do not count this family of solutions here.
For all other values we apply Lemma~\ref{lem:pell-solutions}, which gives
$O_\ep(N^\ep)$ solutions for each ratio. Therefore, in total, we obtain
$O_\ep(N^\ep)$ solutions.
\end{proof}

\begin{proof}[Proof of Lemma~\ref{lem:rademacher-fourth-moment}]
We consider three cases.

\noindent\textbf{Case 1: trivial solutions.} First, we count the trivial solutions, i.e.\ those in which $a,b,c,d$ split into
pairs and, in each pair, the two numbers are equal. (The case where some pair sums
to $N$ is impossible since all variables are $<\frac{N}{2}$.) There are $3$ possible
pairings, and for each pairing we have $T(N)^2$ solutions.
Therefore, the total number of trivial solutions is $3T(N)^2+O(N)$, where the error term comes from overcounting the solutions where $a=b=c=d$.

    We now show that all remaining solutions contribute
$O_\ep(N^{2-\frac{4}{11}+\ep})$.

\medskip
\noindent\textbf{Case 2: non-trivial solutions with a repeated variable.}
Suppose there is one pair of equal numbers; without loss of generality, $a=b$.
Then the total number of solutions is $O_\ep(N^{1+\ep})$: there are $O(N)$ choices
for the trivial pair $(a,b)$, and for each such choice the pair $(c,d)$ has
$O_\ep(N^\ep)$ possibilities by Lemma~\ref{lem:rademacher-second-moment}.

\medskip
\noindent\textbf{Case 3: $a,b,c,d$ are pairwise distinct.}  Assume without loss of generality that $d$ is the largest of the four variables. This assumption affects the final bound only by a factor of $4$.

Define the \emph{contributions} of $a$, $b$, $c$, $d$ as
\[
A=\frac{a(N-a)}{\gcd(a,N)^2},\quad
B=\frac{b(N-b)}{\gcd(b,N)^2},\quad
C=\frac{c(N-c)}{\gcd(c,N)^2},\quad
D=\frac{d(N-d)}{\gcd(d,N)^2}.
\]
Then $A,B,C,D$ are square-free and $ABCD$ is a square.

    We record a simple consequence of Lemma~\ref{lem:pell-solutions}.
    \begin{lem}\label{lem:radem_two_fixed}
       For fixed $a$ and $d$ there are up to $O_\ep(N^\ep)$ possible values of $(b,c)$ such that $b\neq c$, $b+c\neq N$, and the quadruple $(a,b,c,d)$ satisfies the hypothesis of Lemma~\ref{lem:rademacher-fourth-moment}.
    \end{lem}
    
    \begin{proof}[Proof of the Lemma~\ref{lem:radem_two_fixed}]
        Since $A,B,C,D$ are square-free and $ABCD$ is a square,
\[
\frac{AD}{\gcd(A,D)^2}=\frac{BC}{\gcd(B,C)^2}.
\]

Given $a$ and $d$, the values of $A$ and $D$ are fixed. Moreover,
$\frac{B}{\gcd(B,C)}\mid\frac{AD}{\gcd(A,D)^2}$ and
$\frac{AD}{\gcd(A,D)^2}\leqslant N^4$, so it has $O_\ep(N^\ep)$ divisors. Therefore,
$\frac{B}{\gcd(B,C)}$ can take $\ll_\ep N^\ep$ values, and the same holds for
$\frac{C}{\gcd(B,C)}$. Consequently, $\frac{B}{C}$ can take only $\ll_\ep N^\ep$ values.

Furthermore,
\[
\frac{B}{C}
=\frac{b(N-b)}{c(N-c)}\cdot \frac{\gcd(c,N)^2}{\gcd(b,N)^2}.
\]
Since $N$ has $O_\ep(N^\ep)$ divisors, it follows that $\frac{b(N-b)}{c(N-c)}$ also
takes only $\ll_\ep N^\ep $ values, all satisfying the hypotheses of
Lemma~\ref{lem:pell-solutions}. Therefore, for each ratio different from $1$ we
obtain $O_\ep(N^\ep)$ solutions, yielding $O_\ep(N^\ep)$ possibilities for $(b,c)$
when $b(N-b)\neq c(N-c)$.
    \end{proof}
    We now consider three subcases.
    \begin{enumerate}
        \item If $d\leqslant N^{\frac{9}{11}}$.

Then $d$ has $N^{\frac{9}{11}}$ possibilities, and for a fixed $d$ we have
$a\leqslant d$, hence $a$ also has at most $N^{\frac{9}{11}}$ possibilities. Given
$(a,d)$, Lemma~\ref{lem:radem_two_fixed} gives $O_\ep(N^\ep)$ possibilities for
$(b,c)$. Thus, there are $O_\ep(N^{\frac{18}{11}+\ep})$ solutions in total.
        \item If $d> N^{\frac{9}{11}}$ and $D>N^{\frac{12}{11}}$.

Since $D$ is square-free and $ABCD$ is a square, we have $D\mid ABC$, and hence
\[
\gcd(A,D)\gcd(B,D)\gcd(C,D)\geqslant D.
\]
Therefore, at least one of these gcd's is $\geqslant \sqrt[3]{D}>N^{\frac{4}{11}}$.
Without loss of generality, assume $\gcd(A,D)>N^{\frac{4}{11}}$ and set
$Y=\gcd(A,D)$. Since $Y\mid D$ and $D$ is coprime to $N$, we have $\gcd(Y,N)=1$.
Applying Lemma~\ref{lem:Y-divides-aNa}, we obtain that, for fixed $d$, there are
$O_\ep(N^{\frac{7}{11}+\ep})$ possible values of $a$.

Given $(a,d)$ we again have $O_\ep(N^\ep)$ possibilities for $(b,c)$ by Lemma~\ref{lem:radem_two_fixed}, and there are
$O(N)$ possibilities for $d$. Hence, we obtain $O_\ep(N^{\frac{18}{11}+\ep})$
solutions in this case.
\item If $d>N^{\frac{9}{11}}$ and $D\leqslant N^{\frac{12}{11}}$.

For such $d$ we have $d(N-d)\gg N^{\frac{20}{11}}$. Therefore, if
$D\leqslant N^{\frac{12}{11}}$ then
\[
\gcd(N,d)^2\gg N^{\frac{8}{11}}
\qquad\text{and hence}\qquad
\gcd(N,d)\gg N^{\frac{4}{11}}.
\]
There are $O_\ep(N^\ep)$ possible values for $\gcd(d,N)$, and for each such value
there are $O(N^{\frac{7}{11}})$ possibilities for $d$. Thus, there are
$O_\ep(N^{\frac{7}{11}+\ep})$ choices for $d$.

For each such $d$ there are $O(N)$ possibilities for $a$, and given $(a,d)$ there
are $O_\ep(N^\ep)$ possibilities for $(b,c)$ by Lemma~\ref{lem:radem_two_fixed}.
Therefore, this case also contributes $O_\ep(N^{\frac{18}{11}+\ep})$ solutions.
    \end{enumerate}

    Summing the error terms in all cases yields the desired bound.
\end{proof}
Having established the second- and fourth-moment inputs
(Lemmata~\ref{lem:rademacher-second-moment} and~\ref{lem:rademacher-fourth-moment}),
we now conclude Theorem~\ref{thm:main-rademacher} via McLeish's theorem
(Theorem~\ref{McLeish_CLT}). Let us denote by $\Ppsf(r)$ the largest prime divisor of $r$ which has an odd exponent, and assume that $\Ppsf(r)=0$ if $r$ is a perfect square. As in \cite{chinis2025random}, define
\[
M_N=\sum_{n<\frac{N}{2}}f(n)f(N-n),
\qquad
M_{p,N}=\sum_{\substack{n<\frac{N}{2}\\\Ppsf(n(N-n))=p}}f(n)f(N-n),
\]
and their normalised counterparts
\[
S_N=\frac{1}{\sqrt{\EE[M_N^2]}}M_N,
\qquad
S_{p,N}=\frac{1}{\sqrt{\EE[M_N^2]}}M_{p,N}.
\]

It is easy to see that $(S_{p,N})_{p>0}$ form a martingale difference sequence, and $M_{0,N}=0$. We now verify the three conditions of Theorem~\ref{McLeish_CLT}.

\medskip
\noindent\textbf{Normalised variances.}
We rewrite the condition as follows:
\begin{align*}
\sum_{p\leqslant N}\EE[S_{p,N}^2]
&=\frac{1}{\EE[M_N^2]}\sum_{p\leqslant N}\EE[M_{p,N}^2]\\
&=\frac{1}{\EE[M_N^2]}\sum_{p\leqslant N}
  \sum_{\substack{n_1,n_2<\frac{N}{2}\\\Ppsf(n_1(N-n_1))=\Ppsf(n_2(N-n_2))=p}}
  \EE[f(n_1)f(N-n_1)f(n_2)f(N-n_2)]\\
&=\frac{1}{\EE[M_N^2]}\#
\Bigl\{
\begin{aligned}[t]
&n_1,n_2<\tfrac{N}{2}\ \Bigm|\ 
n_1,n_2,N-n_1,N-n_2\text{ are square-free},\\
& n_1(N-n_1)n_2(N-n_2)=\square,\ 
\Ppsf(n_1(N-n_1))=\Ppsf(n_2(N-n_2))
\Bigr\}.
\end{aligned}
\end{align*}
If $n_1=n_2$ and $n_1,N-n_1$ are square-free, then
$\Ppsf(n_1(N-n_1))=\Ppsf(n_2(N-n_2))$ and such pairs contribute to the count. If not,
then by Lemma~\ref{lem:rademacher-second-moment} there are $O_\ep(N^\ep)$ such pairs
even without the condition on $\Ppsf(\cdot)$, hence also $O_\ep(N^\ep)$ with this
condition. Therefore,
\[
\sum_{p\leqslant N}\EE[S_{p,N}^2]
=\frac{1}{\EE[M_N^2]}\cdot\bra{T(N)+O_\ep(N^\ep)}.
\]
Similarly, $
\EE[M_N^2]=T(N)+O_\ep(N^\ep)$. Since $T(N)\asymp N$, we obtain $
\sum_{p\leqslant N}\EE[S_{p,N}^2]=1+o(1)$,
and the normalised variances condition holds.

\medskip
\noindent\textbf{Fourth-moment condition.}
We rewrite:
\begin{align*}
\sum_{p\leqslant N}\EE[S_{p,N}^4]
&=\frac{1}{\EE[M_N^2]^2}\sum_{p\leqslant N}\EE[M_{p,N}^4]\\
=&\frac{1}{\EE[M_N^2]^2}\sum_{p\leqslant N}
  \sum_{\substack{n_1,n_2,n_3,n_4<\frac{N}{2}\\\Ppsf(n_1(N-n_1))=\dots=\Ppsf(n_4(N-n_4))=p}}
  \EE[f(n_1)f(N-n_1)\dots f(n_4)f(N-n_4)]\\
=&\frac{1}{\EE[M_N^2]^2}\#
\Bigl\{
\begin{aligned}[t]
&n_1,n_2,n_3,n_4<\tfrac{N}{2}\ \Bigm|\ 
n_1,N-n_1,\dots,n_4, N-n_4\text{ are square-free},\\
&\qquad n_1(N-n_1)\dots n_4(N-n_4)=\square,\\
&\qquad \Ppsf(n_1(N-n_1))=\dots =\Ppsf(n_4(N-n_4))
\Bigr\}.
\end{aligned}
\end{align*}
By Lemma~\ref{lem:rademacher-fourth-moment}, there are $o(N^2)$ quadruples
$(n_1,n_2,n_3,n_4)$ satisfying the square-freeness and square conditions and not
splitting into equal pairs. Moreover, we have $\EE[M_N^2]\gg N$. Hence, to prove the
fourth-moment condition, it suffices to verify that
\[
\#\Bigl\{
n,m<\tfrac{N}{2} \Bigm|
n,m,N-n,N-m\text{ are square-free; }
\Ppsf(n(N-n))=\Ppsf(m(N-m))
\Bigr\}
=o(N^2).
\]
In fact, we prove the following stronger statement.
\begin{lem}\label{lem:radem-smooth-aux}
\[
\#\{n,m<\tfrac{N}{2}\mid n,m\text{ are square-free and }\Ppsf(n(N-n))=\Ppsf(m(N-m))\}=o(N^2).
\]
\end{lem}
\begin{proof}
    We use the same trick as in \cite{Klurman2023}, \cite{chinis2025random}, and in the
    proof of Theorem~\ref{thm:main-steinhaus} above.

    Consider two cases.
    \begin{itemize}
    \item $\Ppsf(n(N-n))=\Ppsf(m(N-m))=p<\sqrt{\log N}$.

Then, since $n$ is square-free, $\frac{n}{\gcd(n,N^\infty)}$ is $\sqrt{\log N}$-smooth, and by formula $(1.18)$ from \cite{granville08},
the number of possibilities for $\frac{n}{\gcd(n,N^\infty)}$ is at most
\[
\bra{\frac{\log N}{\sqrt{\log N}}}^{(1+o(1))\pi\bra{\sqrt{\log N}}}.
\]
Moreover,
\begin{align*}
\bra{\frac{\log N}{\sqrt{\log N}}}^{(1+o(1))\pi\bra{\sqrt{\log N}}}
&=\bra{\log N}^{\frac{1+o(1)}{2}\cdot \frac{\sqrt{\log N}}{\frac{1}{2}\log\log N}}\\
&=\exp\bra{(1+o(1))\sqrt{\log N}}
=O_\ep(N^\ep).
\end{align*}
Also, by Claim~\ref{claim:bound-radicals}, there are $\ll_\ep N^\ep$ possibilities for $\gcd(n, N^\infty)$. So, there are $\ll_\ep N^\ep$ possibilities for $n$. Similarly for $m$, and hence there are $O_\ep(N^\ep)=o(N^2)$ pairs in this case.
\item $\Ppsf(n(N-n))=\Ppsf(m(N-m))=p\geqslant\sqrt{\log N}$.

Then either $n$ or $N-n$ is divisible by $p$. For a fixed $p$, there are at most
$\frac{2N}{p}$ possibilities for $n$, and similarly for $m$. Therefore, the total
number of pairs is bounded by
\[
\sum_{p\geqslant\sqrt{\log N}}\frac{4N^2}{p^2}
\leqslant \sum_{n\geqslant\sqrt{\log N}}\frac{4N^2}{n^2}
=O\bra{\frac{N^2}{\sqrt{\log N}}}
=o(N^2).
\]
\end{itemize}

Combining the two cases proves the lemma.
\end{proof}

Applying Lemma~\ref{lem:radem-smooth-aux}, we obtain
\[
\sum_{p\leqslant N}\EE[S_{p,N}^4]=o(1).
\]

\medskip
\noindent\textbf{Cross-term condition.}
As in the proof of Theorem~\ref{thm:main-steinhaus}, once the fourth-moment
condition is established, it is enough to verify that
\[
\EE\sqbra{\bra{\sum_{p\leqslant N}S_{p,N}^2}^2}=1+o(1).
\]
We can rewrite
\begin{align*}
\EE\sqbra{\bra{\sum_{p\leqslant N}S_{p,N}^2}^2}
&=\frac{1}{\EE[M_N^2]^2}\#
\Bigl\{
\begin{aligned}[t]
&n_1,n_2,n_3,n_4<\tfrac{N}{2}\ \Bigm|\ 
n_1,N-n_1,\dots,n_4,N-n_4\text{ are $\square$-free},\\
&\qquad n_1(N-n_1)\dots n_4(N-n_4)=\square,\\
&\qquad \Ppsf(n_1(N-n_1))=\Ppsf(n_2(N-n_2)),\\
&\qquad \Ppsf(n_3(N-n_3))=\Ppsf(n_4(N-n_4))
\Bigr\}.
\end{aligned}
\end{align*}
There are $T(N)^2$ trivial solutions with $n_1=n_2$ and $n_3=n_4$. The other two trivial pairing families give $o(N^2)$, because then $\Ppsf(n_i(N-n_i))$ are the same for all $i$, and we can apply Lemma~\ref{lem:radem-smooth-aux}. By
Lemma~\ref{lem:rademacher-fourth-moment}, there are $o(N^2)$ further solutions (as
before, we drop the conditions involving the greatest prime divisor when counting
non-trivial solutions). Moreover, we showed above that
\[
\EE[M_N^2]=T(N)+O_\ep(N^\ep).
\]
Therefore,
\[
\EE\sqbra{\bra{\sum_{p\leqslant N}S_{p,N}^2}^2}
=\frac{1}{(T(N)+o(N))^2}\cdot\bra{T(N)^2+o(N^2)}.
\]
Since $T(N)\gg N$, it follows that
\[
\EE\sqbra{\bra{\sum_{p\leqslant N}S_{p,N}^2}^2}=1+o(1),
\]
which is the desired cross-term condition.

\medskip
\noindent\textbf{Conclusion.}
Having verified all three conditions of Theorem~\ref{McLeish_CLT}, we can conclude that
$\sum_{p\leqslant N}S_{p,N}$ converges in distribution to $\mathcal{N}(0,1)$.
Therefore,
\[
\frac{1}{\sqrt{\EE[M_N^2]}}\cdot\sum_{n<\frac{N}{2}}f(n)f(N-n)
\xrightarrow[]{d}\mathcal{N}(0,1).
\]
Finally, since $\EE[M_N^2]=T(N)+o(T(N))$, we obtain
\[
\frac{\sum_{n<\frac{N}{2}}f(n)f(N-n)}{\sqrt{T(N)}}
\xrightarrow[]{d}\mathcal{N}(0,1),
\]
which completes the proof of Theorem~\ref{thm:main-rademacher}.
\section*{Acknowledgements}

I would like to thank Dr Oleksiy Klurman for proposing this problem and for his extensive mathematical guidance on this project. I would also like to express my deepest gratitude for his invaluable mentorship and support as I transition to postgraduate research. \blue{I also thank Besfort Shala for his helpful comments on an earlier draft of this paper.}
\printbibliography
\end{document}